\documentclass[twoside,11pt,reqno]{amsart}

\usepackage{upref,amsxtra,amssymb,amscd}
\usepackage{varioref}
\usepackage{verbatim}
\usepackage{epsfig}

\usepackage{epic,eepic,eucal}
\usepackage{enumerate}

\def\B{          \mathcal B}

\def\A{          \mathcal A}
\def\mS{          \mathcal S}
\def\O{          \mathcal O}
\def\Id{        \rm Id}

\def\a{         \alpha}

\def\N{ {\bf N}   }

\def\Z{ {\bf Z}   }

\newcommand{\NN}{{\mathbb N}}
\newcommand{\RR}{{\mathbb R}}
\newcommand{\TT}{{\mathbb T}}
\newcommand{\ZZ}{{\mathbb Z}}
\newcommand{\CC}{{\mathbb C}}
\newcommand{\QQ}{{\mathbb Q}}
\newcommand{\DD}{{\mathbb D}}
\newcommand{\BB}{{\mathbb B}}
\def\carre{ \hfill $\Box$    }

\newtheorem{theo}{\sc Theorem}[section]
\newtheorem{prop}[theo]{\sc Proposition}
\newtheorem{lemm}[theo]{\sc Lemma}
\newtheorem{coro}[theo]{\sc Corollary}

\theoremstyle{definition}

\newtheorem{defi}[theo]{\sc Definition}

\theoremstyle{remark}

\newtheorem{rema}[theo]{\sc Remark}

\newtheorem{prob}[theo]{\sc Problem}

\numberwithin{equation}{section}

\begin{document}
\title[Constructions in elliptic dynamics.]{Constructions in elliptic dynamics.}
\author[Bassam Fayad, Anatole Katok]{Bassam Fayad, Anatole Katok*}
\address{Bassam  Fayad, Universit\'e Paris 13, CNRS 7539, Villetaneuese 93430  }
\email{fayadb@math.univ-paris13.fr}
\address{Anatole Katok, Department of Mathematics,
  Penn State University, State College, PA 16802}
\email{katok\_a@math.psu.edu}
\dedicatory
{Dedicated to the memory of Michel Herman}

\begin{abstract}

We present an overview and  some new applications of the approximation by 
conjugation method  introduced by Anosov and the second author more than 
thirty years ago \cite{AK}.   Michel Herman  made  important contributions  to the development and applications  of this method  beginning from the construction of minimal and uniquely ergodic diffeomorphisms jointly with Fathi in  \cite{FH} and continuing with exotic invariant sets of rational maps of the Riemann sphere \cite{H3}, and the construction of invariant tori with  nonstandard and  unexpected behavior in the context of KAM theory
\cite{H1, H2}. Recently the method  has been experiencing a revival. Some of the new results presented in the paper illustrate 
variety of  uses for  tools available for a long time, others  exploit  new   methods,  in particular possibility of mixing in the context of Liouvillean  dynamics discovered by the first author \cite{F1, F2}.

\end{abstract}

\thanks{*) Partially
supported by NSF grant DMS-0071339}
\maketitle

\section{Elliptic dynamics. Diophantine and Liouvillean behavior}
In this paper  we present and discuss a variety of examples  of 
$C^{\infty}$ volume preserving diffeomorphisms 
of compact manifolds which  by any reasonable count  should be  viewed as  exotic if not pathological. The methods which produce those examples are discussed in the next section and specific constructions  carried out in the rest of the paper.
Before doing that we  would like to present a brief  justification  for paying attention to  these 
admittedly very special situations by putting them into a certain general perspective
along the lines of \cite[Section 7]{HK}, where  a more detailed discussion and appropriate references can  be found.

The term ``elliptic dynamics'' is often, although not universally, used to  denote the study of
recurrent (conservative) behavior in dynamical systems at the other end of the spectrum from ``chaotic''.  While in the everyday usage the opposite of ``chaotic''  would be ``ordered'',  in the dynamics context the preferred counterpart is ``stable''. The stable asymptotic behavior in the context of  conservative  motions is best represented by  quasiperiodic motions (translations and linear flows)   on tori, which appear for example  in completely integrable Hamiltonian  systems. Such a system splits on  most of its  phase space into invariant tori  which are determined by the fixed values of the action integrals in the action--angle coordinates. On each torus the motion  is linear in the angle coordinates with the frequency vector determined by the  action vector (see \cite[Section 1.5]{KH} for an introductory discussion). It is  quite remarkable 
that under a suitable nondegeneracy assumption  many (in fact, most in terms of  phase volume)  invariant tori survive  slightly distorted  under any small Hamiltonian  perturbation of the system.
This is the  principal content of KAM theory  to which Michel Herman made major contributions
\cite{H4}. This prevalence of stability  which  also appears in other  contexts (see e.g. \cite{HIHES})  is of course 
a major manifestation of elliptic phenomena.   

Stability appears when the linear model,
the translation (or rotation if multiplicative notations are used) of the  torus $\TT^k$,
$$
T_{\a}x=x+\a,\,\,\,\,\,\,\text{mod}\,\, 1
$$
or the linear flow on $\TT^k$,
$$
T^t_{\omega}x=x+\omega t ,\,\,\,\,\,\,\text{mod}\,\, 1,\,\,\, t\in\RR
$$
is of a \emph{Diophantine} type, i.e.  the translation vector $\a$ or the frequency vector $\omega$  is not too well approximated by rationals.  A typical Diophantine  condition  in the discrete--time case is
$$
|\langle \a, m\rangle- n| > C\|m\|^{-s}
$$
for any $m\in\ZZ^k,\,\,n\in\ZZ$ and some positive numbers $C$ and $s$ independent of 
$m$ and $n$. Similarly in the continuous--time case  one has  for some  positive $C$ and  $s$
and any $m\in\ZZ^k$,
$$
|\langle \omega, m\rangle| > C\|m\|^{-s}.
$$

The opposite case of  exceptionally good rational approximation of  $\a$ or $\omega$ is usually called  \emph{Liouvillean}. The Liouvillean behavior is associated with instability:
circle diffeomorphisms with Liouvillean rotation numbers   usually   have singular invariant measures and  hence are conjugate to the corresponding rotations  via nondifferentiable continuous maps, time--changes of linear flows of 
tori  with Liouvillean frequency vectors often and in fact typically  exhibit unexpected ergodic
properties \cite{F3}, the invariant tori of completely integrable Hamiltonian  systems with Liouvillean frequency vectors  usually do not survive after perturbations (or  may survive with greatly distorted dynamics, see \cite{H1}). 

This leads to another possible opposite to the chaotic behavior. Since the latter is associated with exponential  growth of orbit complexity with time (positive entropy, exponential growth of the number of periodic orbits, etc; see e.g. \cite[Section 6]{HK})  one may consider  slow growth of various characteristics of such  complexity  as a hallmark of elliptic dynamics.  Notice that    translations and linear
 flows of a torus  are isometries, hence there is no divergence of orbits, and in a  certain sense  no
 growth of orbit complexity with time at all.  In the Diophantine case   this behavior  is often stable
 under small perturbations of the system while in the Liouvillean case   perturbations of the linear models  often produce substantive qualitative effects but only after long periods of time.  To understand what kind of effects may result form this instability   a constructive approach  is quite useful.  While there is no chance to  describe
 all possible   effects which appear  in nonlinear perturbations of  Liouvillean  behavior, the variety of possibilities is astonishing. In    a number of cases  such perturbations provide the only known way to exhibit certain types of dynamical behavior, e.g. ergodic behavior accompanied  by only finitely many periodic points on such simple manifolds as   the two--dimensional disc or two--dimensional sphere. See Section \ref{Discdiffeo}   for a more detailed discussion of those cases.

We would like to thank A. Fathi   for pointing out to the paper \cite{Ha}  
where a construction similar to those of Section 4   is carried out, M. 
Handel for valuable comments and observations which 
helped us to improve presentation in  that section,  A. Windsor and M. 
Saprykina for useful discussions. The second author  would like to thank  H. Nakayma
for useful comments  which helped to improve the statement of  Proposition \ref{propNakayama} and the Japan Society for Promotion of Science for supporting the visit  to Japan during which this paper was completed.

\section{An overview of the approximation by conjugation method }
\subsection{General scheme.}\label{general scheme}
 Let $M$ be a differentiable manifold with a nontrivial  smooth circle action $\mS=\{S_t\}_{t\in\RR},\,\, S_{t+1}=S_t$ preserving a smooth volume.  Every  smooth  $S^1$ action preserves a smooth volume $\nu$ which  can  be obtained  by taking  any volume  $\mu$ and averaging  it with respect  to  the action: $\nu=\int_0^1(S_t )_*\mu dt$.  Similarly $\mS$ preserves  a smooth Riemannian metric on $M$ obtained by averaging  of a  smooth Riemannian 
metric.  

When this  does not cause an ambiguity we will often call   the elements $S_{\a}$ of the circle action $\mS$
"rotations", correspondingly "rational" and "irrational",  for rational or irrational values of $\a$.

Denote by  $C_q$ the subgroup of $S^1$ with $q$ elements, i.e the
$q$th roots of unity.
Assuming without loss of generality that  the action is effective (no 
element acts as identity) we see that  
for an open dense subset of $M$ the stationary subgroup is trivial. 
Using the standard properties of isometries one sees the set $F$ 
of fixed points of the action  is a  smooth submanifold of $M$. The fixed 
point 
set $F_q$  of the map $S_{1\over q}$ (or, equivalently, of the action
 of $C_q$)  is also a smooth  submanifold.

Volume preserving  maps  with various interesting, often surprising, topological and  ergodic properties are obtained as limits  of volume preserving periodic transformations 

\begin{eqnarray} \label{conjugations} f = \lim_{n \rightarrow \infty} f_n,\,\,\,\,\text{where}  \,\,\,\, f_n=H_n S_{\a_{n+1}} H_n^{-1} \end{eqnarray} 
with $\a_n = {p_n \over q_n} \in \QQ$ and

\begin{eqnarray} \label{conjugations2} H_n =  h_1 \circ ... \circ h_n,
\end{eqnarray}
where every $h_n$ is a  volume preserving diffeomorphism of $M$ that  satisfies 
\begin{eqnarray} \label{conjugations3} h_n \circ S_{\a_n} = S_{\a_n} \circ h_n.
\end{eqnarray} 

In certain versions of the method the diffeomorphisms $h_n$  are chosen not preserving  the volume but  distorting it   in a controllable way;  this for example is the only interesting situation when $M$ is  the circle (see e.g. \cite[Section 12.6]{KH}).

Usually at step $n$, the diffeomorphism $h_n$ is  
constructed first,  and  $\a_{n+1}$ is chosen 
afterwards close enough to $\a_n$ to guarantee convergence of 
the construction. For example, it is easy to see that for the 
limit in (\ref{conjugations}) to exist in the $C^\infty$ topology it is  largely
sufficient to ask  that 
\begin{eqnarray} \label{convergence} 
|\a_{n+1} - \a_{n}| \leq\frac {1}{ 2^n q_n{||H_n||}_{C^n}}.
\end{eqnarray} 

The power and fruitfulness  of the method   depend on the fact that the sequence of  
diffeomorphisms $f_n$ is made to converge while the conjugacies  $H_n$ diverge
often ``wildly'' albeit in a controlled (or prescribed)  way.
Dynamics of the circle actions and of their individual elements is  simple and well--understood. In particular, no element of such an action  is  ergodic or topologically transitive, unless the circle action itself is transitive, i.e  $M=S^1$.  To provide interesting  asymptotic properties 
of  the limit  typically the successive conjugacies  spread the orbits of the circle action 
$\mS$ (and 
hence  also those of its restriction to the subgroup $C_q$ for any sufficiently large $q$)  across the phase space $M$
making them  almost dense, or almost uniformly distributed, or approximate another type of interesting asymptotic behavior. Due to the high speed of convergence this remains true for sufficiently long orbit segments of the limit diffeomorphism. To guarantee an appropriate speed of approximation  extra conditions on convergence of approximations in addition to (\ref{convergence}) may be required. 

There are many variations of the construction within this general scheme.  In the subsequent sections we describe several representative  versions, each addressing a problem   which has not been treated  before. Presenting each version  to the last detail and in the greatest generality would be  quite tedious. Accordingly we  often choose to  describe routine steps only briefly
and sometimes  treat representative special cases indicating  generalizations  afterwards.
  
\subsection{Almost all vs. all orbits}\label{almostallall}
In different versions of the  approximation by conjugation method one may control 
the asymptotic behavior of  almost all  orbits  with respect to the invariant volume, or  of all orbits.
Somewhat imprecisely we will call those versions ergodic and topological.

Ergodic constructions  deal with   measure--theoretic (ergodic) properties 
with respect to a given invariant volume, such as the  number of ergodic components 
(in particular ergodicity),
rigidity,
weak mixing, mixing, further spectral properties. Topological constructions deal 
with minimality, number of ergodic invariant measures (e.g. unique ergodicity)
and their supports, presence of particular invariant sets, and so on.
The paper \cite{AK} dealt only with ergodic constructions. 
A topological version  was announced in  \cite{K2} (see also \cite{Br})
   and independently developed in \cite{FH}. 

Control over behavior of the orbits of approximating periodic diffeomorphisms
$f_n$  in (\ref{conjugations})  on the $n$th step of the construction is typically provided by taking an
invariant under $S_{\a_n}$ (and hence under $S_{1\over{q_n}})$
 collection of "kernels",
usually smooth balls,  and  redistributing  them  in the phase space in  a prescribed fashion (also $S_{1\over{q_n}}$ invariant).  In  ergodic constructions one requires   the complement to the union  of the kernels  to have small volume and hence most orbits of $\mS$ (and consequently  of  any  finite subgroup $C_q$  for a sufficiently large $q$) to spend most of the time inside the kernels. In the topological versions the kernels need to be chosen in such a way that \emph{every} orbit of $\mS$ spends most of the time inside the kernels.  This   requires more care and  certain attention to the geometry of orbits. See for example Figures 3.1 and 3.2.

A natural way  of selecting the  kernels, their intended images,  and constructing  a map $h_n$ satisfying (\ref{conjugations3})
is by taking  a fundamental domain  $\Delta$ for $S_{\a_n}$ (or, equivalently, for $S_{1\over{q_n}})$
choosing kernels and images inside  $\Delta$, constructing  a diffeomorphism of $\Delta$ to itself 
identical near its boundary  which sends  kernels into their intended images, and extending the map to the images $S_{k\over{q_n}},\,\,\,k=1,\dots,q_n-1$ by  commutativity.
This  method in particular is used in  the construction of a diffeomorphism conjugate to a rotation, see \ref{nonstandardrotation} below, as well as later in this paper in the proofs of Theorems \ref{theoremthreemeasures}, \ref{theoremmain} and \ref{yytheorem}. However in certain cases, for example to  achieve weak mixing, it is necessary to use more general constructions.

\subsection{Closure of conjugates}\label{closurespace} By  controlling initial steps of the construction one can  keep  the  final diffeomorphism $f$ within any given
$C^{\infty}$ neighborhood of  a given element $S_{\a}$ of the action $\mS$.
On the other hand,  by applying  a fixed  diffeomorphism $g$  first  the construction can be carried out  in a neighborhood of  any diffeomorphism  conjugate to an element of the action.  Thus the
closure $\A$ of  diffeomorphisms of the form $gS_tg^{-1}$ in, say, $C^{\infty}$ topology is a natural
ambient  space for  considering what a construction  of this type can produce.
This was first noticed in \cite[Section 7]{AK} in connection with ergodic properties
with respect to the invariant volume and  was used in \cite{FH} to control topological properties.

If the action  $\mS$ has fixed points it is  natural to consider a smaller  space  $\A_F\subset \A$, namely, the closure of   $gS_tg^{-1}$ where
the conjugacy $g$ fixes every point of the fixed point set $F$ of the action $\mS$.
Obviously for $f\in\A_F$ every point of $F$ is also fixed. There does not seem to be much 
difference between dynamical properties which may appear in the whole space $\A$
and in $\A_F$.  In Section
\ref{Discdiffeo} we will discuss the  situation for  the standard two--dimensional 
examples: the  disc, the annulus and the sphere with the $S^1$ action by rotations .  At this point we would like to introduce also \emph{restricted} spaces of conjugacy closures,  namely  for  $\a\in\RR/\ZZ$ we define the spaces
$\A_{\a}$  and $\A_{F,\a}$ as the closure in $C^{\infty}$ topology of the diffeomorphisms of the form $gS_{\a}g^{-1}$ where $g$ is  correspondingly an arbitrary $C^{\infty}$ diffeomorphism, or a $C^{\infty}$ diffeomorphism fixing every point of the set $F$.
Since for any ${p\over q}\in\QQ$ and any $f\in A_{p\over q}$ one has $f^q=\Id$, only spaces $A_{\a}$ for irrational $\a$'s are of interest.

Let us describe a classical example  showing that the space $\A$ may be quite different from the union of $\A_{\a}$

\subsubsection{The horocycle flow}\label{horocycleflow} Consider the group $SL(2,\RR)$ and    three one--parameter subgroups: the standard rotations $R_t=\begin{pmatrix}
\cos 2\pi t & & \sin 2\pi t  \\  -\sin 2\pi t & & \cos 2\pi t \end{pmatrix}$,  the diagonals
$G_t=\begin{pmatrix}  e^{t} & & 0 \\ 0 & & e^{-t }\end{pmatrix}$
 and the upper--triangular  unipotents $H_t=\begin{pmatrix}  1 & & t \\ 0 & & 1 \end{pmatrix}$.  The commutation relation
$$ G_tR_sG_{-t}=\begin{pmatrix}
\cos 2\pi s & & e^{2t}\sin 2\pi s  \\  -e^{-2t}\sin 2\pi s & & \cos 2\pi s \end{pmatrix}
$$
implies
\begin{eqnarray}\label{eqcommutation}
\lim_{t\to\infty}G_tR_{e^{-2t}\over 2\pi}G_{-t}=H_1.
\end{eqnarray}
Now let $\Gamma$ be a cocompact lattice in $SL(2,\RR)$.  Since left translations on any group commute with right translations we can consider the  actions of the  one--parameter subgroups described above on the compact manifold $M=SL(2,\RR)/\Gamma$.
In the case when $\Gamma$ has no elliptic elements, $M$ can be identified with 
the unit tangent bundle of   the compact  surface $V=SO(2)\backslash SL(2,\RR)/\Gamma$ of constant negative curvature, and the right actions of $\{ G_t\}$ and $\{H_t\}$ are called the \emph{geodesic flow} and \emph{horocycle flow} on $V$ correspondingly. \cite[Sections 5.4 and 17.5]{KH}.
We will denote the actions of $R_t,\,\,G_t$ and $H_t$  by right translations on
$SL(2,\RR)/\Gamma$ by $r_t,\,\,g_t$ and $h_t$ correspondingly. The first of these is of course an $S^1$ action, while every nontrivial  element  of the last one is uniquely ergodic, mixing of all orders and has countable Lebesgue spectrum (see \cite[Section 1.4e]{KSS} and references thereof).
It follows from (\ref{eqcommutation}) that 
$$
\lim_{t\to\infty}g_tr_{e^{-2t}\over 2\pi}g_{-t}=h_1.
$$
Thus $h_1\in\A$; it is the limit of conjugacies of  certain elements of the circle action
which converge to zero. 

\begin{prop}
$h_1\notin\A_{\a}$ for any  $\a\in S^1$.
\end{prop}
\begin{proof}
Obviously $\A_0=\{\Id\}$; thus, we need to prove  the statement  for 
$\a\neq 0$;
For any volume preserving diffeomorphism $f$ of a compact manifold  homotopic to identity one defines the \emph{rotation class} $\rho(f)$ as an element of the first homology  group of the direct product 
$M\times S^1$. It is constructed as the asymptotic cycle of the suspension flow of $f$ \cite[Section 14.7b]{KH}.  The rotation class is  continuous  in $C^0$ topology and equivariant under the conjugacy. Thus $\rho(h_1)=0$.  On the other hand the projection of $\rho(r_\a)$ to  $H_1(M,\RR)$ has the form $\a \rho_0$, where $\rho_0$ is the nonzero  integer homology class    realized by  the embedding  of  $SO(2)$ to $SL(2,\RR)/\Gamma$.
Since the action of the diffeomorphism group by conjugacies  respects the integer structure in this projection
the rotation classes of elements of  conjugates of $r_{\a}$ are separated from zero and hence
for any $f\in\A_{\a},\,\,\rho(f)\neq 0$. 
\end{proof}

The following question remains open:
\begin{prob}\label{mixinginA}
Does there exist a circle action $\mS$  on a  compact manifold   and $\a$  such that the  corresponding space $\A_{\a}$ contains a mixing  transformation?
\end{prob}

In Section \ref{sectionmixing} we will give  a positive answer  for a version  of this question  with a $\TT^3$ action instead of an $S^1$ action. We will also provide more examples  of circle actions with mixing 
transformations  in the unrestricted space $\A$.

\subsection{Generic and non-generic constructions}
The development of the approximation by  conjugation method was  motivated 
by several natural questions  about ergodic properties  of volume preserving  diffeomorphisms.
Here are two  characteristic examples  of such questions, both answered  positively in \cite{AK}:

\subsubsection{Does there exist   an  ergodic area preserving diffeomorphism of the two--disc $\DD^2$?}\label{discergodic} $ $

\subsubsection{Does there exist a diffeomorphism of  a  compact manifold $M, \,\, \dim M>1$ 
preserving a smooth volume and measurably conjugate to an irrational rotation of the circle?}\label{nonstandardrotation} $ $
\bigskip

These positive answers  looked somewhat surprising at the beginning: in the former case because  usually near an elliptic  fixed point there are  many invariant curves making the map   highly non-ergodic;  and in the latter because the  eigenfunctions of any such transformation must be  highly discontinuous.
While   ergodic diffeomorphisms  of the disc  without  genuine elliptic points which exhibit much stronger  stochastic properties than ergodicity were found later  \cite{K}, the original examples in \cite{AK} do have elliptic points and are in fact small perturbations of (non-ergodic) rotations of the disc.

Constructions answering  questions  \ref{discergodic} and \ref{nonstandardrotation}  are characteristic for two different   versions of the approximation by conjugation method.  

\subsubsection{Control at some time scales} In the constructions of the first kind it is sufficient  \emph{to control the behavior 
of approximating and hence resulting diffeomorphisms on  a series  of growing but unrelated time scales}. The approximate  pictures may look quite whimsical (see e.g. the original weak mixing construction in \cite[Section 5]{AK} and a modern version in \cite{GK}),
but as long as a diffeomorphism $f\in\A$ is close enough to conjugates of  rotations appearing in such pictures the property is guaranteed. In fact, for a proof of genericity 
in $\A$  of a property exhibited by a construction of this sort no actual inductive construction is needed. One just needs to show that an 
approximate picture  at each scale appears for an open dense  subset of conjugates of rotations. If  appearance in an approximate picture at infinitely many growing scales guarantees the property then by the Baire category theorem  the property holds for
a dense  $G_{\delta}$ subset on $\A$. 

Ergodicity is generic in $\A$; in the case of the  disc with $S^1$ acting by rotations around the center this  in particular implies  that   invariant  circles  which are very hard to avoid from the classical standpoint  (twists, frequency locking, etc) can  in fact be destroyed  easily if one adopts  a less conventional point of view.  We will discuss  the disc diffeomorphisms in detail is Section
\ref{Discdiffeo} below.

Beyond ergodicity the basic ergodic properties which are generic in $\A$  are  arbitrary fast  cyclic approximation (which implies  ergodicity and also  simple singular spectrum) 
\cite[Theorem 7.1]{AK}, and weak mixing. Those  properties are  among   generic 
 properties  for  measure preserving automorphisms of Lebesgue space in weak 
topology \cite{KS, K1} which also  possess  some more refined  generic 
properties such 
as mutual singularity of the maximal spectral type and all its convolutions 
\cite{K1}, 
or existence of roots of all orders \cite{Ki}.  The first of those  as well 
as many 
others follows from the genericity of  arbitrary fast periodic approximation of   arbitrary type (not just cyclic)  \cite{K1}. 

We will later discuss open questions  and difficulties related with
 establishing genericity of such properties in $\A$ (Section \ref{sectionsmooth}).

In this paper we present  a new example of  a generic construction in 
Section \ref{sectionthreemeasures}. A proper extension of this construction  
from the basic  two--dimensional examples (the disc, the annulus and the sphere)
   to  manifolds (with or without boundaries) with circle actions provides   a  
common generalization of the  genericity of ergodicity  in the general  case 
of nontrivial circle actions \cite{AK}, and the genericity of minimality  and 
unique ergodicity in the case of fixed--point free actions on closed 
manifolds \cite{FH} (See Section \ref{minimalgeneral}).

\subsubsection{Control at all times} In the constructions of the second kind  including the one answering the question     
\ref{nonstandardrotation} approximations  at different  steps of the construction are linked and hence in principle \emph{the asymptotic behavior of the resulting diffeomorphism is controlled for  all times}. Such a construction is never generic in $\A$ or any of the other spaces introduced in 
Section \ref{closurespace}

Constructions of this kind appear most naturally  when  the resulting 
diffeomorphism  is constructed to be measure--theoretically conjugate to a  
map of a particular kind such as a rotation in 
the question \ref{nonstandardrotation} or a  translation on a torus of higher 
dimension
\cite[section 6]{AK}. Those constructions  also appear in  topological 
situations, for 
example in producing   minimal diffeomorphisms   with more than one ergodic measure \cite{W}.
Our construction in Section \ref{Minimalsets}  is non-generic although quite natural in the sense described at the end of Section \ref{almostallall},
and is  similar in some respects to   the construction of a diffeomorphism conjugate to a rotation.
One may also extend the  construction from \cite{W} to obtain extra invariant 
measures for the case
of  actions  with fixed points or  on manifolds with boundaries (such as  
four or more ergodic measures on the disc; see Section 
\ref{sbsmoremeasures}). The 
construction in 
Section \ref{sectiontrnonerg} is also non-generic and is of a similar type.

A new and quite different type of non-generic constructions are those  which 
produce mixing diffeomorphisms. Presently we  need an action of $\TT^3$ to 
carry  out such a construction.
This construction is   based on the  combination of  approximation by 
conjugation and  the methods  developed by the first author for constructing  
mixing time changes for linear flows on $\TT^k,\,\,k\ge 3$ \cite{F1, F2}. The latter may be interpreted as a ``multiple frequencies''
version of the  approximation by conjugation method.

We would like to point out that Michel Herman was  a great master in the use of Baire category  theorem  with results which were often striking (it was his "favorite theorem" according to Albert Fathi), and in particular he very effectively used 
 generic constructions in the context of the method of approximation by conjugation.
In his published work though he  did not make any direct use of non-generic constructions
in the context of the approximation by conjugation method.
It was however his original insight which stimulated the first author, his last Ph. D. student,
to  develop the non-generic constructions of mixing diffeomeorphisms.

\section{Diffeomorphisms of the disk, the annulus, and the sphere  with exactly 
three ergodic measures.}\label{sectionthreemeasures} 
\subsection{Dynamics and invariant measures of disc 
diffeomorphisms}\label{Discdiffeo}
The  most basic fact about  dynamics of diffeomorphisms of the disc $\DD^2$
is presence  of a fixed point (the Brouwer Fixed point theorem).   Furthermore,  by the Brouwer Translation  theorem  for an  area preserving
map there is fixed point inside the disc. 

As early as 1930  L.G. Shnirelman \cite{Sh} in an attempt to prove a 
"converse"  to the Brouwer theorem 
constructed an example of  a topologically transitive homeomorphism $h$ of 
$\DD^2$ which in  somewhat modified polar coordinates  $(\rho,\theta)$ has  the form
\begin{eqnarray} \label{discskewproduct} 
h(\rho,\theta)=(\rho+g(\theta), \theta +\a).
\end{eqnarray}
(Here, for example one can take $\rho=-(\tan(r/\pi))^{-1}$, where $r$ is the usual  radius).
Shnirelman was mistaken in his assertion that  every orbit inside the disc other than the origin was dense. For maps of the form (\ref{discskewproduct})  existence of nondense orbits in the interior  can be easily shown. (A.S. Besicovitch who repeated this mistake later corrected it in \cite{Be}). It is also not too difficult to see 
that for any homeomorphism of the disc fixing the origin there are  nondense \emph{semiorbits} 
in the interior. Existence of nondense \emph{orbits} is much more difficult and 
it was proven only recently (forty five years after the Besicovitch paper) by P. Le Calvez and J.-C. Yoccoz \cite{CY}.

It is true however that the only ergodic  invariant measures in the Shnirelman example 
are
the $\delta$--measure at the origin and Lebesgue measure on the boundary and  
that  for  almost every  point with respect to the two--dimensional Lebesgue measure   the asymptotic distribution  exists and  is  exactly the  average of these two measures.
A careful study of the Shnirelman--Besicovitch examples was one of the sources which
led Anosov and the second  author  to the  invention of the approximation by conjugation method (\cite[Introduction]{AK}).

Notice that  unlike Shnirelman--Besicovitch maps any area preserving  diffeomorphism of the disc  has at least three ergodic invariant measures: the $\delta$--measure at a fixed point inside the disc, a measure  supported at the boundary plus any ergodic component  of the area.
Naturally, the only  possibility to achieve this lower bound of three ergodic invariant measures appears
when there is only one fixed point inside, only one measure supported  by the boundary and
the map is ergodic with respect to the area (this is  not sufficient though).

Now consider the  closure spaces  $\A$  and $\A_F$ of conjugates of rotations  described in 
Section \ref{closurespace}.  Every  diffeomorphism  $f\in\A_F$ preserves the origin, the boundary and the volume. 
In this section we will  show that the minimal number of three ergodic invariant measures  is achieved and  is even  generic in a properly defined sense. Moreover in our examples (albeit not generically) 
any point in the interior whose orbit closure is not dense (such points exist by 
\cite{CY})  has an asymptotic distribution which is either the $\delta$ measure at the origin,
or Lebesgue measure on the boundary, in contrast with the Shnirelman--Besicovitch
maps of the form (\ref{discskewproduct}).  

Any map $f\in\A$ preserves the boundary and hence has a rotation number on the 
boundary.
By continuity for $f\in\A_{\a}$ this rotation number is equal to $\a$. 
Our examples belong to spaces $\A_{\a}$ with   $\a$ irrational but very well  
approximated by rationals (Liouvillean numbers). This is the only possibility 
for ergodicity with respect to the area
in a space $\A_{\a}$. For, if $\alpha$ is rational then every $f\in\A_{\a}$ 
is periodic.  If  $\a$ is Diophantine the restriction of every $f\in\A_{\a}$ 
to the boundary of $\DD^2$  is $C^{\infty}$
conjugate  to the rotation $R_{\a}$ by the Herman--Yoccoz theorem \cite{HIHES, Y}.
Hence as Herman proved in his unpublished work there are smooth $f$--invariant circles 
arbitrary close to the boundary and hence ergodicity or topological transitivity is not possible.


The only other way ergodicity with  respect to the area may appear in the whole space $\A$  would be for  a map with a rational rotation number
$p\over q$ on the boundary which is a limit of  conjugacies of rotations  with variable rotation numbers converging to $p\over q$.  While we do not know whether 
such ergodic examples exist their structure must be quite different since due to recent results of 
Le Calvez \cite{LC2003}  any such map has 
infinitely many periodic points (and hence infinitely many ergodic invariant measures) in the interior. \footnote{This has been pointed out to us by the referee of the paper}

The following problem is related to the above discussion as well as to  Problem \ref{mixinginA} .
\begin{prob}
Does there exist an area preserving  mixing diffeomorphism of the disc $\DD^2$ with zero metric (or topological) entropy?
\end{prob}

\begin{rema}Notice that such examples  exist  on  any   orientable surface of genus $g\ge 1$  with or without holes.  They are obtained as time one maps of  mixing flows first constructed by Ko\v cergin
\cite{kochergin}. It is not  known  whether  there are  examples on the sphere or the cylinder.
\end{rema}

\subsection{Statement and description of the construction}\label{theorem 
three measures}
Notice that the pictures of rotations and conjugacies  are essentially 
identical on the disc, 
the annulus $[0,1]\times S^1$ and the sphere $S^2$. In all three cases  one  
has polar coordinates $(\rho, \theta)$ and the rotations are given by 
$S_t(\rho,\theta)=(\rho,\theta +t)$.
Here $0\le \rho \le 1$ and $\theta$ is  the cyclic coordinate of period 1. 
All the conjugacies  involved in our constructions will  coincide with the 
identity near $r=0$ and $r=1$, hence the difference between the three cases will be insignificant. Accordingly our drawings will
 refer to the case of the annulus  which is easier to represent graphically.
Let us call the Lebesgue measure  on the manifold, the $\delta$--measures
at the fixed points of the rotations and Lebesgue measures on the boundary 
components the  \emph{natural measures}.

\begin{theo}\label{theoremthreemeasures}Let $M$ be $\DD^2,\,\,[0,1] \times 
S^1 $ or $S^2$, and $S_t$ be the standard action by  rotations.
There exists a $C^\infty$ diffeomorphism $f: M\to M$   that has exactly 
three ergodic invariant measures, namely  the natural measures on $M$. 
Furthermore, diffeomorphisms with this  property  form a residual set  in the
space $\A'$: the closure in the  $C^{\infty}$ topology of the  conjugates of 
rotations
with conjugacies fixing the fixed points of $\mS$  and  every point of the 
boundary.
\end{theo}

We will call the boundary components and fixed points of the rotations  \emph{singularities}.
At step $n-1$, $h_{n-1}$ is constructed and then $\a_n$ is chosen. 
We assume that after step $n-1$ there are two neighborhoods of the singularities $B_{n-1}^l$ and $B_{n-1}^r$ on which $H_{n-1} \equiv$ Id. 
At step $n$ we define three collections of closed sets $\delta^{l}_n=\{ \delta^{l}_{n,i}, i=1,\dots, l_n\},\,\, \delta^{l}_{n,i} \subset B^l_{n-1}$, $\delta^{r}_n=\{ \delta^{r}_{n,i}, i=1,\dots, l_n\},\,\, \delta^{r}_{n,i}\subset B_{n-1}^r$ and $\xi_n=\{ \xi_{n,i}, i=1,\dots, l_n\}$, the collection of kernels. These collections will have the property that the union of their atoms capture an increasingly large proportion of \emph{any} rotation orbit. The map $h_n$ is then chosen so that the images of  the  kernels  $H_n (\xi_{n,i})$ have small diameters,  and such that $H_n \equiv $ Id on two smaller neighborhoods  of the singularities $B_{n}^l$ and $B_{n}^r$ that still contain all elements of  $\delta^l_n$ and $\delta^r_n$ correspondingly. At the end of step $n$, $\a_{n+1}$ is chosen to insure convergence of the construction. We will use polar coordinates $(\rho, \theta)$ in our calculations.

\subsection{ Properties of the collections $\xi_n$, $\delta^{r}_n$ and $\delta^r_n$.} \label{collections}

The following proposition  gives a   description of the configuration 
represented   by Figure 3.1.
Notice that the parameter   $\eta_n$ describing precision of the configuration 
can be made arbitrarily small since we are free to take $l_n$ arbitrarily 
large. 

\begin{prop} \label{collections step n} Assuming that neighborhoods  
$B^l_{n-1}$ and $B^r_{n-1}$ of the boundary components are  given and 
given any $\eta_n > 0$, there exists 
an integer $l_n$ multiple of $q_n$ and three collections of disjoint closed sets on $M$, $\delta^l_n$, $\delta^r_n$ and $\xi_n$ each one containing $l_n$ elements and satisfying 

\begin{enumerate}  

\item \label{collections 1}  $\delta^l_n$, $\delta^r_n$ and $\xi_n$ are invariant by $S_{1 \over l_n}$ and $\xi_{n,1} \subset (0,1) \times (0, {1 \over l_n})$,  $\delta^l_{n,1} \subset [0, 1 / n] 
\times (0, {1 \over l_n})$ and   $\delta^r_{n,1} \subset [1- 1 / n,1] \times (0, {1 \over l_n})$,

\item \label{collections 2} $\mu (\xi_{n,1}) \leq \eta_n $ and $\mu(\bigcup \xi_{n,i}) \geq 1 - \eta_n$.

\item \label{collections 3}  $\delta^l_{n,i} \subset B^{l}_{n-1}$ and $\delta^r_{n,i} \subset B^{r}_{n-1}$

\item \label{collections 4} Any circle $V_{\rho} = \{ \rho \} \times S^1$ intersects the three collections in the following way: there exists a constant $c_n( \rho) \in [0,1]$ with $c_n ( 1 / 2) =1$, $c_n(0)=c_n(1)=0$ and such that for any $1 \leq i \leq l_n$,  $V_{\rho} \bigcap \xi_{n,i}$ is an empty set, a point or a segment and so are $ V_\rho \bigcap \delta^l_{n,i}$ and $V_\rho \bigcap \delta^r_{n,i}$; moreover we require that

\begin{eqnarray*} \label{xi} (1 - \eta_n) {c_n(\rho) \over l_n}  \leq |V_\rho \bigcap \xi_{n,i}| \leq  (1 + \eta_n) {c_n( \rho ) \over l_n}, \end{eqnarray*}
and if $\rho \in [0, 1 / 2]$, $V_\rho \cap \delta^r_n = \emptyset$ while  

\begin{eqnarray*} \label{deltal} (1 - \eta_n) {(1 - c_n( \rho )) \over l_n}  
\leq |V_\rho \bigcap \delta^{l}_{n,i}| \leq  (1 + \eta_n) {(1- c_n( \rho)) 
\over l_n},  \end{eqnarray*}
and if  $\rho \in [1 / 2, 1]$, $V_\rho \cap \delta^l_n = \emptyset$ while  
 
\begin{eqnarray*} (1 - \eta_n) {(1 - c_n( \rho )) \over l_n}  \leq |V_\rho 
\cap \delta^r_{n,i}| \leq  (1 + \eta_n) {(1- c_n( \rho )) \over l_n}. 
\end{eqnarray*}

\end{enumerate}

\end{prop}

\vspace{1cm}

\begin{center}

\setlength{\unitlength}{0.00041667in}
\begingroup\makeatletter\ifx\SetFigFont\undefined%
\gdef\SetFigFont#1#2#3#4#5{%
  \reset@font\fontsize{#1}{#2pt}%
  \fontfamily{#3}\fontseries{#4}\fontshape{#5}%
  \selectfont}%
\fi\endgroup%
{\renewcommand{\dashlinestretch}{30}
\begin{picture}(6975,4206)(0,-10)
\path(1275,4158)(1275,33)(6075,33)(6075,4158)
\path(6075,108)(5475,108)(5400,108)(5625,483)
\path(5625,483)(6075,483)
\path(1275,108)(1725,108)(1950,483)(1275,483)
\path(2025,483)(1800,108)(5325,108)
	(5550,483)(2025,483)
\dashline{60.000}(1275,558)(6075,558)
\path(1275,633)(1725,633)(1950,1008)(1275,1008)
\path(2025,1008)(1800,633)(5325,633)
	(5550,1008)(2025,1008)
\dashline{60.000}(1275,1083)(6075,1083)
\path(6075,633)(5475,633)(5400,633)(5625,1008)
\path(5625,1008)(6075,1008)
\dashline{60.000}(1275,2733)(6075,2733)
\path(1275,2283)(1725,2283)(1950,2658)(1275,2658)
\path(2025,2658)(1800,2283)(5325,2283)
	(5550,2658)(2025,2658)
\path(6075,2283)(5475,2283)(5400,2283)(5625,2658)
\path(375,33)(375,2208)
\path(1050,33)(1050,558)
\path(300,2208)(450,2208)
\path(300,33)(450,33)
\path(975,33)(1125,33)
\path(975,558)(1125,558)
\path(5625,2658)(6075,2658)
\thicklines
\path(1275,4158)(1275,33)
\path(6075,4158)(6075,33)
\thinlines
\path(1275,2208)(6075,2208)
\dashline{60.000}(5025,33)(5025,35)(5025,39)
	(5024,47)(5024,59)(5023,76)
	(5022,99)(5021,128)(5020,163)
	(5018,203)(5016,249)(5014,299)
	(5012,353)(5010,411)(5008,470)
	(5005,532)(5003,594)(5001,656)
	(4999,717)(4997,778)(4995,837)
	(4993,894)(4992,950)(4990,1004)
	(4989,1055)(4988,1105)(4987,1153)
	(4986,1199)(4986,1243)(4985,1286)
	(4985,1327)(4985,1368)(4985,1407)
	(4985,1446)(4985,1484)(4986,1521)
	(4987,1558)(4988,1595)(4989,1635)
	(4990,1674)(4991,1713)(4993,1753)
	(4994,1793)(4996,1833)(4998,1874)
	(5000,1916)(5002,1957)(5004,2000)
	(5007,2042)(5009,2085)(5012,2128)
	(5014,2172)(5017,2215)(5020,2259)
	(5022,2302)(5025,2346)(5028,2389)
	(5030,2432)(5033,2475)(5036,2517)
	(5038,2558)(5041,2599)(5043,2639)
	(5046,2679)(5048,2718)(5050,2756)
	(5052,2793)(5054,2829)(5056,2865)
	(5057,2900)(5059,2934)(5060,2967)
	(5061,3000)(5063,3033)(5063,3069)
	(5064,3105)(5065,3141)(5065,3177)
	(5065,3214)(5065,3250)(5065,3288)
	(5064,3327)(5063,3366)(5062,3408)
	(5061,3450)(5059,3495)(5058,3541)
	(5056,3588)(5054,3637)(5051,3687)
	(5049,3738)(5046,3790)(5043,3840)
	(5041,3890)(5038,3937)(5036,3981)
	(5033,4021)(5031,4056)(5029,4086)
	(5028,4111)(5027,4129)(5026,4143)
	(5025,4151)(5025,4156)(5025,4158)
\dashline{60.000}(2295,32)(2295,34)(2295,37)
	(2295,44)(2294,55)(2294,70)
	(2293,91)(2293,118)(2292,150)
	(2291,188)(2289,233)(2288,283)
	(2287,339)(2285,399)(2283,464)
	(2282,533)(2280,604)(2278,678)
	(2276,754)(2274,830)(2272,907)
	(2271,983)(2269,1059)(2267,1133)
	(2266,1206)(2264,1277)(2263,1347)
	(2261,1414)(2260,1479)(2259,1542)
	(2258,1602)(2257,1661)(2257,1717)
	(2256,1772)(2256,1824)(2255,1875)
	(2255,1925)(2255,1972)(2255,2019)
	(2255,2064)(2255,2108)(2255,2152)
	(2256,2195)(2256,2237)(2257,2278)
	(2258,2319)(2258,2364)(2259,2409)
	(2261,2453)(2262,2498)(2263,2542)
	(2265,2587)(2267,2632)(2268,2678)
	(2270,2725)(2273,2772)(2275,2821)
	(2278,2871)(2281,2922)(2284,2974)
	(2287,3029)(2291,3084)(2294,3142)
	(2298,3201)(2302,3261)(2307,3323)
	(2311,3386)(2316,3449)(2320,3513)
	(2325,3577)(2330,3640)(2334,3702)
	(2339,3762)(2343,3819)(2348,3873)
	(2351,3923)(2355,3969)(2358,4010)
	(2361,4046)(2363,4076)(2365,4101)
	(2367,4121)(2368,4135)(2369,4146)
	(2370,4152)(2370,4155)(2370,4157)
\put(0,1083){\makebox(0,0)[lb]{\smash{{{\SetFigFont{8}{9.6}{\familydefault}{\mddefault}{\updefault}${1 \over q_n}$}}}}}
\put(5325,3558){\makebox(0,0)[lb]{\smash{{{\SetFigFont{8}{9.6}{\familydefault}{\mddefault}{\updefault}$B_{n-1}^r$}}}}}
\put(1500,3558){\makebox(0,0)[lb]{\smash{{{\SetFigFont{8}{9.6}{\familydefault}{\mddefault}{\updefault}$B_{n-1}^l$}}}}}
\put(5625,258){\makebox(0,0)[lb]{\smash{{{\SetFigFont{6}{7.2}{\familydefault}{\mddefault}{\updefault}$\delta_{n,1}^r$}}}}}
\put(3225,258){\makebox(0,0)[lb]{\smash{{{\SetFigFont{6}{7.2}{\familydefault}{\mddefault}{\updefault}$\xi_{n,1}$}}}}}
\put(675,183){\makebox(0,0)[lb]{\smash{{{\SetFigFont{8}{9.6}{\familydefault}{\mddefault}{\updefault}${1 \over l_n}$}}}}}
\put(1350,258){\makebox(0,0)[lb]{\smash{{{\SetFigFont{6}{7.2}{\familydefault}{\mddefault}{\updefault}$\delta_{n,1}^l$}}}}}
\put(6975,783){\makebox(0,0)[lb]{\smash{{{\SetFigFont{8}{9.6}{\familydefault}{\mddefault}{\updefault}.}}}}}
\end{picture}
}

\end{center}

$$ {\rm Fig. \ 3.1. \ The \ collections \ } \xi_n, \delta^l_n, {\rm \ and  \ 
} \delta^r_n. $$

\subsection{Properties of the successive conjugating diffeomorphisms $H_n$.}\label{Section3.4} Using the collections of sets introduced in \S \ref{collections}, we will state and prove some sufficient conditions on the sequences $\a_n$ and $H_n$ that guarantee  the convergence in (\ref{conjugations}) to a diffeomorphism having exactly three invariant ergodic measures as in Theorem \ref{theorem three measures}. We require that the collections in Proposition \ref{collections step n} be constructed with $\eta_n \rightarrow 0$ and we will extract later more information on how fast $\eta_n$ should converge to 0.  

\begin{lemm} \label{Hn} If $\a_n = {p_n \over q_n} $ is a sequence of rationals and $H_n$ is a sequence of diffeomorphisms of $M$ satisfying (\ref{conjugations2}) and (\ref{conjugations3}), and if there exists a sequence $\epsilon_n \rightarrow 0$ such that for every $n$ large enough we have

\begin{enumerate}
\item \label{Hn a}  $H_n \equiv {\rm Id}$ on two open neighborhoods $B_n^l$ of $ \lbrace 0 \rbrace \times S^1$ and $B_n^r$ of $ \lbrace 1 \rbrace \times S^1$
satisfying 
$$ \delta^j_{n,i} \subset B^j_n \subset B^j_{n-1}, \ \ j= r, l,\,\,\,i=1,\dots,l_n.$$

\item \label{Hn b} For any $i=1,\dots,l_n$, 
$$ {\rm diam} (H_n (\xi_{n,i})) \leq \epsilon_n,$$ 

\item \label{Hn c} $|\a_{n+1} - \a_n| \leq \frac{1}{ 2^n q_n {||H_n||}^n_{C^n}}$, and $q_{n+1} \geq nl_n$,

\end{enumerate}

\noindent then $H_n \circ S_{\a_{n+1}} \circ H_n^{-1} $ converges in the $C^\infty$ topology to a diffeomorphism $f$ having exactly three invariant ergodic measures: the Lebesgue measure $\mu$ on $M$ and the  measures $\delta^l$ and $\delta^r$  supported on the singularities  corresponding to Lebesgue  measures  on $\lbrace 0 \rbrace \times S^1$ and $\lbrace 1 \rbrace \times S^1$ respectively.

\end{lemm}

\begin{rema} \label{feasability} In \S \ref{const} we will have to construct 
$h_n$ so that $H_n$ fulfills conditions (\ref{Hn a}) and (\ref{Hn b}) above. 
Condition (\ref{Hn a}) is consistent with (\ref{collections 3}) in 
Proposition \ref{collections step n} and condition (\ref{Hn b}) above will be 
consistent with (\ref{collections 2}) in Proposition \ref{collections step n} 
if at step $n$ we choose $\eta_n \leq {(\epsilon_n / {||H_{n-1}||}_{C^1})}^3$.
\end{rema}

\subsubsection{Convergence.}  Condition (\ref{Hn c}) of  Lemma \ref{Hn} 
insures convergence of $f_n = H_n \circ S_{\a_{n+1}} \circ H_n^{-1}$  in the $C^{\infty}$ topology to a volume preserving diffeomorphism $f$ and the choice of $\a_{n+1}$ is made only 
after $h_n$ is constructed. It also insures the closeness between $f^m$ 
and $f_n^m$ as $m \leq q_{n+1}$ (Cf. \S \ref{fn}).  The last condition, 
$q_{n+1} \geq n l_n$, guarantees that for any point in $M$ a $(1 - {1 \over 
\eta_n})-$proportion  of its $q_{n+1}$ (vertical) orbit under 
$S_{\a_{n+1}}$ 
is captured by the atoms of $\xi_n$, $\delta^l_n$ and $\delta^r_n$ which is crucial in the estimations of ergodic averages (Cf. (\ref{22}) and (\ref{33})).

\vspace{0.4cm}

\begin{center}

\setlength{\unitlength}{0.00041667in}
\begingroup\makeatletter\ifx\SetFigFont\undefined%
\gdef\SetFigFont#1#2#3#4#5{%
  \reset@font\fontsize{#1}{#2pt}%
  \fontfamily{#3}\fontseries{#4}\fontshape{#5}%
  \selectfont}%
\fi\endgroup%
{\renewcommand{\dashlinestretch}{30}
\begin{picture}(4920,4098)(0,-10)
\path(33,258)(1308,258)(1758,1008)(33,1008)
\path(33,1233)(1308,1233)(1758,1983)(33,1983)
\path(33,2208)(1308,2208)(1758,2958)(33,2958)
\path(3633,258)(1683,258)(1533,258)
	(1983,1008)(3858,1008)
\path(3633,1233)(1683,1233)(1533,1233)
	(1983,1983)(3858,1983)
\path(3633,2208)(1683,2208)(1533,2208)
	(1983,2958)(3858,2958)
\dashline{60.000}(3033,33)(3033,3858)
\thicklines
\path(33,3708)(33,33)
\path(33,2958)(33,3858)
\thinlines
\drawline(3783,1983)(3783,1983)
\drawline(3783,1008)(3783,1008)
\drawline(3633,258)(3633,258)
\drawline(3483,1233)(3483,1233)
\drawline(3558,1233)(3558,1233)
\drawline(3633,1233)(3633,1233)
\drawline(3633,2208)(3633,2208)
\path(3558,2208)(4683,2208)
\path(3633,1233)(4683,1233)
\path(3633,258)(4683,258)
\drawline(3783,2958)(3783,2958)
\drawline(3783,1983)(3783,1983)
\drawline(3858,1008)(3858,1008)
\path(3783,2958)(4908,2958)
\path(3783,1983)(4908,1983)
\path(3783,1008)(4908,1008)
\thicklines
\path(1308,33)(1309,35)(1313,41)
	(1319,51)(1328,67)(1340,88)
	(1356,115)(1375,148)(1396,185)
	(1420,226)(1445,270)(1471,315)
	(1496,359)(1522,404)(1546,447)
	(1570,487)(1592,526)(1612,562)
	(1631,596)(1649,627)(1665,656)
	(1679,683)(1693,708)(1705,731)
	(1717,752)(1727,772)(1737,790)
	(1746,808)(1760,837)(1772,864)
	(1783,888)(1793,912)(1802,933)
	(1809,953)(1815,971)(1819,987)
	(1822,1002)(1825,1014)(1826,1025)
	(1826,1034)(1825,1042)(1824,1048)
	(1823,1053)(1821,1058)(1817,1063)
	(1814,1068)(1809,1072)(1803,1075)
	(1795,1078)(1787,1080)(1777,1081)
	(1766,1082)(1753,1083)(1739,1083)
	(1724,1083)(1708,1083)(1695,1083)
	(1680,1083)(1664,1083)(1647,1083)
	(1629,1083)(1611,1084)(1592,1084)
	(1572,1085)(1553,1085)(1534,1086)
	(1516,1087)(1499,1089)(1484,1090)
	(1470,1092)(1457,1093)(1446,1095)
	(1438,1097)(1431,1099)(1425,1101)
	(1420,1104)(1415,1107)(1411,1110)
	(1407,1115)(1405,1120)(1403,1125)
	(1402,1132)(1402,1140)(1403,1148)
	(1405,1157)(1407,1167)(1411,1179)
	(1415,1191)(1421,1204)(1427,1217)
	(1433,1232)(1441,1248)(1449,1265)
	(1458,1283)(1467,1301)(1477,1320)
	(1487,1340)(1498,1362)(1510,1384)
	(1523,1408)(1536,1433)(1549,1460)
	(1564,1487)(1578,1514)(1593,1542)
	(1608,1571)(1623,1599)(1638,1627)
	(1652,1654)(1667,1681)(1680,1708)
	(1693,1733)(1706,1757)(1718,1779)
	(1729,1801)(1739,1821)(1749,1840)
	(1758,1858)(1767,1876)(1775,1893)
	(1783,1909)(1789,1924)(1795,1938)
	(1801,1951)(1805,1963)(1809,1974)
	(1811,1985)(1813,1994)(1814,2003)
	(1814,2011)(1813,2018)(1811,2025)
	(1809,2031)(1805,2036)(1801,2040)
	(1796,2044)(1791,2048)(1785,2052)
	(1778,2055)(1771,2058)(1759,2062)
	(1746,2067)(1732,2071)(1717,2076)
	(1700,2081)(1682,2085)(1664,2090)
	(1646,2095)(1627,2101)(1609,2106)
	(1591,2110)(1574,2115)(1559,2120)
	(1545,2124)(1532,2129)(1521,2133)
	(1513,2136)(1506,2139)(1500,2143)
	(1495,2147)(1490,2151)(1486,2156)
	(1483,2162)(1481,2169)(1480,2176)
	(1480,2185)(1481,2195)(1482,2206)
	(1485,2218)(1490,2231)(1495,2245)
	(1501,2261)(1508,2278)(1516,2296)
	(1525,2316)(1535,2337)(1546,2359)
	(1558,2383)(1567,2401)(1576,2419)
	(1586,2439)(1597,2461)(1609,2483)
	(1622,2508)(1636,2534)(1650,2563)
	(1666,2593)(1684,2626)(1702,2661)
	(1722,2699)(1743,2739)(1766,2782)
	(1789,2826)(1814,2873)(1839,2921)
	(1865,2969)(1891,3018)(1916,3066)
	(1941,3112)(1964,3155)(1985,3195)
	(2003,3230)(2019,3260)(2032,3285)
	(2043,3304)(2050,3317)(2054,3326)
	(2057,3331)(2058,3333)
\put(483,1533){\makebox(0,0)[lb]{\smash{{{\SetFigFont{8}{9.6}{\familydefault}{\mddefault}{\updefault}$\delta_n^l$}}}}}
\put(558,3183){\makebox(0,0)[lb]{\smash{{{\SetFigFont{8}{9.6}{\familydefault}{\mddefault}{\updefault}${\bf B_n^l}$}}}}}
\put(1083,3858){\makebox(0,0)[lb]{\smash{{{\SetFigFont{8}{9.6}{\familydefault}{\mddefault}{\updefault}$B_{n-1}^l$}}}}}
\put(3858,1533){\makebox(0,0)[lb]{\smash{{{\SetFigFont{8}{9.6}{\familydefault}{\mddefault}{\updefault}$\xi_n$}}}}}
\end{picture}
}

\end{center}

$$ {\rm Fig. \ 3.2. \ } B_n^l \subset B_{n-1}^l.$$
 
\vspace{0.2cm}

\noindent {\sc Proof of Lemma \ref{Hn}.} For $\varphi  \in C^{\infty} (M, \CC)$ we use the notations

$$
\hat{\varphi} = \int_M \varphi(u,v) du dv, \,\,\,\,\,
\hat{\varphi}^l = \int_{S^1} \varphi(0,v) dv, \,\,\,\,\,
\hat{\varphi}^r = \int_{S^1} \varphi(1,v) dv.
$$

\subsubsection{Criterion for the existence of exactly three ergodic invariant 
measures.} Fix $\varphi \in C^{\infty} (M, \CC)$ and $z \in M$ and denote 
by $S_{f,m} \varphi(z)$ the Birkhoff sums 
$$S_{f,m} \varphi(z)= \varphi(z) + \varphi(f(z))+ ... + \varphi(f^{m-1} (z)).$$ Lemma \ref{Hn} clearly follows if  we prove that there exists a sequence of functions $\varepsilon_n \rightarrow 0$, a sequence $m_n \rightarrow \infty$ and a sequence $c_n(z) \in [0,1]$ such that for either $j = l$ or $j=r$  we have

\begin{eqnarray} \label{main}   c_n(z) \hat{\varphi} + (1-c_n(z)) 
\hat{\varphi}^j - \varepsilon_n  \leq {S_{f, m_n} \varphi (z) \over m_n}   \leq 
 c_n(z)  \hat{\varphi} + (1-c_n(z)) \hat{\varphi}^j + \varepsilon_n.
\end{eqnarray}

\vspace{0.4cm}

\setlength{\unitlength}{0.00041667in}
\begingroup\makeatletter\ifx\SetFigFont\undefined%
\gdef\SetFigFont#1#2#3#4#5{%
  \reset@font\fontsize{#1}{#2pt}%
  \fontfamily{#3}\fontseries{#4}\fontshape{#5}%
  \selectfont}%
\fi\endgroup%
{\renewcommand{\dashlinestretch}{30}
\begin{picture}(13887,3934)(0,-10)
\put(9371.651,-215.656){\arc{11129.373}{3.2221}{3.9176}}
\put(-1095.742,4590.654){\arc{13584.690}{0.1453}{0.7112}}
\put(5550,270){\ellipse{3450}{526}}
\thicklines
\path(3975,307)(5475,3532)(7050,307)(3975,307)
\thinlines
\drawline(13875,3907)(13875,3907)
\path(5400,3682)(5401,3683)(5407,3689)
	(5419,3700)(5433,3713)(5446,3724)
	(5457,3733)(5466,3740)(5475,3744)
	(5482,3748)(5489,3750)(5497,3752)
	(5504,3753)(5513,3754)(5521,3753)
	(5528,3752)(5536,3750)(5543,3748)
	(5550,3744)(5559,3740)(5568,3734)
	(5577,3726)(5586,3718)(5594,3709)
	(5602,3700)(5608,3691)(5613,3682)
	(5616,3673)(5619,3664)(5621,3653)
	(5623,3640)(5624,3626)(5625,3614)
	(5625,3608)(5625,3607)
\path(3825,232)(3825,231)(3825,225)
	(3826,213)(3827,199)(3829,186)
	(3831,175)(3834,166)(3838,157)
	(3842,148)(3848,139)(3856,130)
	(3864,121)(3873,113)(3882,105)
	(3891,99)(3900,94)(3907,91)
	(3914,89)(3922,87)(3929,86)
	(3938,85)(3946,86)(3953,87)
	(3961,89)(3968,91)(3975,94)
	(3984,99)(3993,106)(4004,115)
	(4017,126)(4031,139)(4043,150)
	(4049,156)(4050,157)
\put(7425,232){\makebox(0,0)[lb]{\smash{{{\SetFigFont{9}{10.8}{\familydefault}{\mddefault}{\updefault}$\delta^r$}}}}}
\put(5850,3382){\makebox(0,0)[lb]{\smash{{{\SetFigFont{9}{10.8}{\familydefault}{\mddefault}{\updefault}$\delta^l$}}}}}
\put(3225,232){\makebox(0,0)[lb]{\smash{{{\SetFigFont{9}{10.8}{\familydefault}{\mddefault}{\updefault}$\mu$}}}}}
\put(0,2032){\makebox(0,0)[lb]{\smash{{{\SetFigFont{5}{6.0}{\familydefault}{\mddefault}{\updefault}.}}}}}
\end{picture}
}

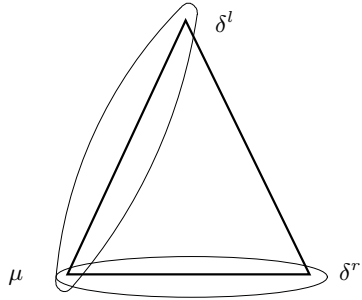
$$ {\rm Fig. \ 3.3.   \ The \ simplex \ of \ invariant \  measures \ and \ 
possible \ averages \ \ for \ a \ large \  } m_n. $$

\vspace{0.2cm}

\subsubsection{The diffeomorphisms $f$ and $f_n = H_n \circ S_{\a_{n+1}} \circ H_n^{-1} $.} \label{fn}  It follows from condition (\ref{Hn c}) of Lemma \ref{Hn} that 

\begin{lemm} \label{fnn} For any $m \leq q_{n+1}$ we have 
$$||f^m - f^m_n||_{C^1} \leq {1 \over 2^n}.$$ 
\end{lemm}

\noindent {\sc Proof.} Since $h_{n+1}$ commutes with $S_{\a_{n+1}}$ we have that $f_n =  H_{n+1} \circ S_{\a_{n+1}}\circ H_{n+1}^{-1} $ hence 
\begin{eqnarray*} 
||f_{n+1}^m - f_n^m||_{C^1} \leq {||H_{n+1}||}_{C^1} {||H_{n+1}||}_{C^2}
m |\a_{n+2} - \a_{n+1}| \leq  \frac{m} { 2^{n+1} q_{n+1}}
\end{eqnarray*} 
by (3) of Lemma 3.5 which yields the conclusion of the lemma.
 \carre

As a consequence (\ref{main}) will follow (with some $\varepsilon_n' \rightarrow 0$) if we pick $m_n = q_{n+1}$ and prove that for a sequence $\varepsilon_n \rightarrow 0$ we have

\begin{eqnarray} \label{main2}   c_n \hat{\varphi} + (1-c_n) \hat{\varphi}^j  - \varepsilon_n \leq {S_{f_n, q_{n+1}} \varphi (z) \over q_{n+1}} \leq 
 c_n \hat{\varphi} + (1-c_n) \hat{\varphi}^j + \varepsilon_n.
\end{eqnarray}

\subsubsection{Control of a large proportion of any orbit.} The point $z \in M$ being fixed we write it as $H_n(\rho_0, \theta_0)$ so that

$$f_n^m (z) = H_n \circ S_{\a_{n+1}}^m (\rho_0, \theta_0).$$
For a set $A \subset M$ we define 
$$L_{q_{n+1}}(A) :=  \lbrace 0 \leq m \leq q_{n+1} -1  \ / \ S_{\a_{n+1}}^m (\rho_0, \theta_0) \in A \rbrace.$$

The following lemma is a straightforward consequence of (\ref{collections 4}) of Proposition \ref{collections step n}

\begin{lemm} \label{distribution} We have for every $1 \leq i \leq l_n$,

 \begin{eqnarray} \label{number 1} (1 - \eta_n) {q_{n+1} c_n( \rho_0) \over   l_n} - 1 \leq \# L_{q_{n+1}} (\xi_{n,i}) \linebreak \leq (1 + \eta_n) {q_{n+1} c_n( \rho_0) \over   l_n} + 1,  \end{eqnarray} 
 and
 \begin{eqnarray}
 \label{number 2}\begin{aligned} (1 - \eta_n){ q_{n+1} (1-c_n( \rho_0))  \over l_n} - 1 \leq \# L_{q_{n+1}} (\delta^j_{n,i}) \\   \leq (1 + \eta_n) { q_{n+1} (1-c_n( \rho_0))  \over l_n}+ 1, 
 \end{aligned} \end{eqnarray} 
where $j=l$ if $\rho_0 \leq {1 / 2}$ and $j=r$ if $\rho_0 \geq {1 / 2}$.

  In particular, the number of points among the first $q_{n+1}$ iterates 
of $(\rho_0, \theta_0)$ under $S_{\a_{n+1}}$ that are not considered in the 
lemma is less than $2 \eta_n q_{n+1} +4l_n$.
  \end{lemm}

\subsubsection{ Estimating $\varphi$ on $H_n (\xi_n)$.} Proposition 
\ref{collections step n} and the fact that $S_t$ and $H_n$ preserve the area 
 $\mu$ imply that $ 1 - \eta_n \leq l_n \mu(H_n(\xi_{n,i})) < 1$ for any $i \leq l_n$.  Hence, as a consequence of condition (\ref{Hn b}) in Lemma \ref{Hn} there exists a sequence $\epsilon_{1,n} \rightarrow 0$ such that for any $(\rho, \theta) \in H_n(\xi_{n,i})$ we have

\begin{eqnarray} 
(1 - \epsilon_{1,n} ) \int_{H_n(\xi_{n,i})}  \varphi(u,v) du dv \leq {1\over l_n} \varphi(\rho, \theta) \nonumber \\ \label{2kernel} 
 \leq  
(1 + \epsilon_{1,n} ) \int_{H_n(\xi_{n,i})}  \varphi(u,v) dudv.   \end{eqnarray}

On the other hand property (\ref{collections 2}) in 
Proposition \ref{collections step n} and the fact that $H_n$ preserves $\mu$  
also implies that there exists $\epsilon_{2,n} \rightarrow 0$ such that

\begin{eqnarray} \label{22kernel} (1 - \epsilon_{2,n} ) \hat{\varphi} \leq
\sum_{i=1}^{l_n} \int_{H_n(\xi_{n,i})}  \varphi(u,v) dudv  \leq  (1 + \epsilon_{2,n} ) \hat{\varphi}, 
\end{eqnarray}
where we recall that $\hat{\varphi}= \int_{M} \varphi(u,v) du dv$.

\subsubsection{Estimating $\varphi$ on $H_n (\delta_n^l)$ and $H_n (\delta_n^r)$.} As a consequence of (\ref{collections 1}) in Proposition \ref{collections step n} , there exists a sequence $\epsilon_{3,n}$ such that for any $(\rho, \theta) \in \delta^{l}_{n,i}$ we have that 
\begin{eqnarray} (1 - \epsilon_{3,n}) \varphi( 0, {i \over l_n}) \leq \varphi(\rho, \theta) \leq 
 (1 + \epsilon_{3,n}) \varphi(0, {i \over l_n}). \label{2boundary} \end{eqnarray}

On the other hand there exists a sequence $\epsilon_{4,n} \rightarrow 0$ such that 
\begin{eqnarray} \label{22boundary} 
 (1 - \epsilon_{4,n} ) \hat{\varphi}^l  \leq {1 \over l_n} \sum_{i=1}^{l_n} \varphi(0, {i \over l_n}) \leq 
 (1 + \epsilon_{4,n} ) \hat{\varphi}^l,
\end{eqnarray}
where as before  $ \hat{\varphi}^l = \int_{S^1} \varphi(0,v) dv $. Similar equations as (\ref{2boundary}) and   (\ref{22boundary}) hold for $(\rho, \theta) \in  \delta^{r}_{n,i}$ and $ \hat{\varphi}^r = \int_{S^1} \varphi(1,v) dv $.

\vspace{0.2cm}

\subsubsection{Proof of (\ref{main2}).} We will assume for definiteness 
that $z = H_n(\rho_0, \theta_0)$ with $\rho_0 \leq {1 / 2}$ the case 
$\rho_0 \geq 1 / 2$ being similar.
In light of Lemma \ref{distribution} and (\ref{2kernel})--(\ref{22boundary}) we get the following estimations of 
$$U_n := {1 \over q_{n+1}} S_{f_n,q_{n+1}} \varphi (\rho_0, \theta_0) $$

\begin{eqnarray}      V_n^1 + V_n^2 + V_n^3 \leq U_n \leq U_n^1 + U_n^2 + U_n^3
\end{eqnarray}
where $U_n^3$ is a bound on the contribution of the non-controlled points 
\begin{eqnarray} \label{11} U_n^3 \leq  2 (\epsilon_n + {2l_n \over 
q_{n+1}}) {||\varphi||}, \end{eqnarray}
and $U_n^1$ is relative to the points that fall in $H_n(\xi_n)$ for which 
we use 
(\ref{number 1}) of Lemma \ref{distribution} and (\ref{2kernel}) and 
(\ref{22kernel}) to obtain  the following upper bound on their contribution 
\begin{eqnarray}  U_n^1 = {1 \over q_{n+1} } \sum_{i=1}^{l_n}  \left[ (1 + 
\eta_n){c_n( \rho_0) q_{n+1} } +  l_n \right] (1 + \epsilon_{1,n}
) \int_{H_n(\xi_{n,i})}  \varphi(u,v) dudv \nonumber \\ \label{22}
\leq  (1 + \eta_n) (1 + \epsilon_{1,n} ) (1 + \epsilon_{2,n} ) c_n(
\rho_0)  \hat{\varphi} + {l_n \over q_{n+1}} ||\varphi||,  \ \ \ \ \
\end{eqnarray}

and $U_n^2$ is a bound on the contribution of the points in $\delta_n^l$ for which  (\ref{number 2}) of Lemma \ref{distribution} and (\ref{2boundary}) and (\ref{22boundary}) lead to  
\begin{eqnarray}  U_n^2 =  {1 \over q_{n+1} }  \sum_{i=1}^{l_n} \left[ (1 + \eta_n) {(1-c_n(\rho_0)) q_{n+1} \over l_n} + 1 \right]  (1 + \epsilon_{3,n}) \varphi(0, {i \over l_n}) 
\nonumber \\ \label{33} \leq  (1 + \eta_n) (1 + \epsilon_{3,n})(1 + \epsilon_{4,n} )(1-c_n(\rho_0)) \hat{\varphi}^l + {l_n \over q_{n+1}}  ||\varphi||. \end{eqnarray}
Recall that $\a_{n+1}$ is chosen such that $l_n / q_{n+1} \leq 1 / n$, hence (\ref{11})--(\ref{33}) with similar lower bounds for  $V_n^1$, $V_n^2$ and $V_n^3$ yield (\ref{main2}). This finishes the proof of  Lemma \ref{Hn}. \carre

\subsection{ Construction of $h_n$ and choice of $\a_{n+1}$.} \label{const} In this section we  construct $h_n$ and choose $\a_{n+1}$ so that the conditions of Lemma \ref{Hn} hold for $H_n = h_1 \circ ... \circ h_n$. 

We fix a sequence $\epsilon_n \rightarrow 0$. We assume that the conditions of Lemma    \ref{Hn} hold up to step $n-1$ and consider a configuration of collections as in Proposition \ref{collections step n} with  
\begin{eqnarray} \eta_n = {(\epsilon_n / {||H_{n-1}||}_{C^1})}^3. \label{eeeta} \end{eqnarray}
Then we have

\begin{prop} \label{hn}
 There exists a diffeomorphism $h_n$ of class $C^\infty$ on $M$ and two neighborhoods $B_n^l$ and $B_n^r$ of the singularities with $\delta_n^j  \subset B_n^j \subset B_{n-1}^j$, $j=l, r$, such that 

\begin{enumerate}

\item \label{100} $h_n$ is volume preserving,

\item \label{200} $h_n \circ S_{\a_n} = S_{\a_n} \circ h_n$,

\item \label{300} $h_n \equiv {\rm Id}$ on $B_n^l \bigcup B_n^r$
\item \label{400} For every $\xi_{n,i} \in \xi_n$, ${\rm diam} 
(h_n (\xi_{n,i})) \leq {\epsilon_n \over {||H_{n-1}||}_{C^1}}$.

\end{enumerate}

\end{prop}

\begin{proof} As explained in Remark  \ref{feasability}, (\ref{300}) above is possible due to   (\ref{collections 3}) in Proposition \ref{collections step n} and the fact that the  unions of the sets in the collections $\delta^j_n$, $j= l, r$ are  closed while $B^j_{n-1}$, $j=l, r$ are open. Condition (\ref{200}) can be realized in particular with $h_n \equiv {\rm Id }$ on the singularities of the fundamental domain $[0,1] \times [0, 1 / q_n]$ due to  (\ref{collections 1}) in Proposition \ref{collections step n} and to the fact that $l_n$ is a multiple of $q_n$. Condition (\ref{400}) can be realized since $\mu(\xi_{n,i}) \leq \eta_n$ and that we chose in (\ref{eeeta}) $\eta_n = {(\epsilon_n / {||H_{n-1}||}_{C^1})}^3 = o( {(\epsilon_n / {||H_{n-1}||}_{C^1})}^2)$.  
\end{proof}

\vspace{0.2cm}

\noindent {\sl Checking the conditions of Lemma \ref{Hn}.}  Clearly, Conditions (\ref{Hn a}) and (\ref{Hn b}) on $H_n$ in Lemma \ref{Hn} follow respectively from (\ref{300}) and (\ref{400}) of Proposition \ref{hn}. Finally the choice of $\a_{n+1}$ is done so that condition (\ref{Hn c}) of Lemma \ref{Hn} holds. The construction is thus completed. \carre
\bigskip   

\noindent {\sl Proof of genericity of diffeomorphisms with three ergodic measures in $\A'$}.
Fix  a function $\varphi$ and  $\epsilon>0$, and define  for  an integer $m$ and a function $c$ 
the set $S(\varphi,\epsilon,m,c)$ of maps $f\in\A'$ such that the
averages  at time $m$ are  close to $c \mu + (1-c) \delta$.  
This set is open by definition.
We actually proved that the union over $c$ and  over $m$
is dense hence the intersection over a countable base $\varphi_i$
and a sequence $\epsilon_j\to 0$ is a dense $G _{\delta}$. \carre

\subsection{Minimal number of invariant measures in $\A$}\label{minimalgeneral}
Consider  a general setting for approximation  by conjugation construction, namely a
smooth action  $\mS$ of the circle on a compact manifold $M$, possibly with boundary.
For the special case of fixed point free actions on closed manifolds Fathi and Herman showed \cite{FH} that  there is a residual subset of  $\A$ which consists of uniquely ergodic diffeomorphisms.

A proper  fairly straightforward generalization of the construction  described above  produces in the case of a general  $S^1$ action on a  closed manifold a residual subset of  $\A_F$  which consists of  diffeomorphisms
with only one invariant measure (the smooth volume) not supported on the set $F$.
In other words, ergodic  invariant  measures for such a diffeomorphism are $\delta$-measures at  the points of $F$ and the smooth volume.

The only essential extra observation concerns the  structure of orbits of a 
circle action near 
the set $F$. In order to obtain a proper counterpart of the picture presented 
on Figure 3.2  one needs to take an $\mS$ invariant  neighborhood of a 
connected component of $F$ in place of a thin vertical strip and perturb it in a $S_{1\over q}$ invariant way  to make it  mostly transversal to the orbits.

The general case of a manifold with boundary in dimension greater than two is slightly more complicated. Unlike the disc and the annulus where the action $\mS$ on the  boundary components is transitive
and hence has ergodic elements, in this case 
in order to minimize the number of ergodic measures one needs to carry out the 
construction to the boundary. This can 
be done with a  certain care, the result being a diffeomorphism in $\A_F$ whose only nonatomic ergodic measures are the smooth volume and smooth measures on  connected components of the boundary.
A good test case to  understand this picture is the three--dimensional  closed ball $\DD^3$
with the circle action by rotations around a fixed axis. A perturbation we are referring to 
is  ergodic volume preserving,  ergodic area preserving  on the boundary,  fixes every point of the rotation axis,
and has no other  ergodic invariant measures.  Another  example with a finite number of invariant
measures is a rotationally symmetric  solid  torus (``the doughnut''). Here a small perturbation of a rotation may produce a volume preserving diffeomorphism   with only two ergodic measures,
the volume and a smooth measure on the boundary.

 \subsection{More measures}\label{sbsmoremeasures} One may ask  
how   the number  of 
ergodic invariant measures can be increased in a controlled way. We  will  briefly discuss the  standard cases  of area preserving diffeomorphisms of  the annulus, the disc, and the sphere. Higher dimensional cases are actually easier to handle.
 
 An obvious way to produce an example with any  odd number $2m+1, m\ge 2$ of invariant measures is by
 pasting $m$ annuli,  or  $m-1$ annuli capped on one side  by a disc,  or $m-2$ annuli capped on each side  by a disc,  with  diffeomorphisms described  earlier in this section. One  picks the same rotation number in each copy and  since  our examples are $C^{\infty}$ tangent to the corresponding rotation on the boundary the resulting diffeomorphisms will be $C^{\infty}$.  Naturally these examples are neither  ergodic  nor topologically transitive.
 
 A modification  of our method allows to produce  in the spaces $\A'$  
topologically transitive 
 diffeomorphisms with  exactly  $m\ge 4$  ergodic invariant measures all but 
two of which are absolutely continuous.  These diffeomorphisms are naturally 
not ergodic with  respect to the area measure. Construction remains the 
same  near the  
singularities  while inside instead of a single 
 parallelogram $\xi_{n,1}$ and its shifts one takes $m-2$ next to each 
other  in such a way 
 that their union still intersects every vertical line by a set of large 
conditional measure 
 to guarantee  a counterpart of condition (4) of  Proposition \ref{collections 
step n}. The images of those parallelograms have small diameters (Cf. 
Lemma \ref{Hn}, (3)).  At this step however there is a new element in the 
construction which makes it non-generic. There are $m-2$ invariant cylinders 
 whose images will represent  distributions which we intend to have as 
approximating  ergodic absolutely continuous invariant measures for the 
limit diffeomorphism.  These distributions must converge at least weakly.
 We may in fact  make their densities converge in $L_1$ by keeping most of each cylinder inside itself. However extra steps are needed to provide topological transitivity. This is done  by sending a small number of parallelograms in  each of $m-2$ collection to other cylinders. A  construction of this kind in a more general and considerably more complicated setting is described in detail in \cite{W}.

Finally, it is possible to  increase the  number of ergodic invariant 
measures  in a controlled way  and still keep the volume ergodic; furthermore, the 
additional  singular ergodic measures may be  made supported  by the whole 
manifold.   Corresponding construction involves conjugacies which are not volume preserving  but which at the end produce a volume preserving  diffeomorphism  For reasons of space we will not even outline  this construction here.

\section{ Minimal sets with arbitrary measure.}\label{Minimalsets} 

\subsection{Preliminaries}
In this section  we describe an  answer  to a  question  posed in 1977 in the paper by Fathi and Herman
\cite[Problem 1.6]{FH}.  We do not claim  any particular novelty or originality; rather our goal is to put  the question into the general context of the  approximation by conjugation method and indicate a  variety of possible topological and ergodic properties for  exotic minimal sets. A similar construction  for the disc case can be found in Handel's paper \cite{Ha}. Handel shows that his minimal set is nowhere  locally connected (in fact,  a pseudo-circle)  although he does not address the  question about its  Lebesgue measure. Herman himself produced similar examples for rational maps of the Riemann sphere \cite{H3}.

\begin{theo}\label{theoremmain} Let $M$ be an $m$-dimensional   differentiable manifold $m\ge 2$, $\mu$ a normalized smooth volume on $M$.  Given any number $s\in[0,1)$
there exists a  compactly supported $C^{\infty}$ diffeomorphism $f$ of $M$ preserving $\mu$
and a compact invariant  minimal set $C$ of $f$ such that 
\begin{enumerate}
\item $\mu(C)=s$\label{thmm1}
\item the set $C$ is connected, has dimension $m-1$\label{thmm2}
and is not homeomorphic locally to a product  of $\RR^{m-1}$ with a Cantor set.
\end{enumerate}
Furthermore, $f$ can be chosen arbitrary close to identity  in the $C^{\infty}$ topology.
\end{theo}

As in the previous section we use  an appropriately adjusted version of the approximation  by conjugation method  of  Section \ref{general scheme}. 
By specifying the parameters of the construction extra information  can be provided both about topology of the set $C$ and about topological  and measure--theoretic (ergodic) properties of $f$
restricted to that set. We will  discuss this after  the proof of the theorem in its basic form.

\begin{rema}
If the manifold $M$ admits a nontrivial  smooth action of the circle
a proper modification of our construction   allows to construct  a diffeomorphism with a minimal set with desired properties in any $C^{\infty}$ neighborhood of any element of the 
circle action.
\end{rema}
\subsection{Reduction to $\DD^{m-1}\times S^1$.}\label{reduction}

Consider the product  $P=\BB^{m-1}\times S^1$ of the $m-1$-dimensional  open unit  ball with the circle.
Let  $\lambda$ be the product of normalized  Lebesgue measures. 
We will call  the $m-1$--dimensional direction  ``horizontal''  and  draw it accordingly, and the $S^1$ direction  `` vertical''.  The action of the circle  on $P$   by vertical translations  will  be denoted $\mS=\{S_t\}$ as usual. In  the natural  system of coordinates $(\rho, \theta)$ one
 obviously   has  $S^t(\rho, \theta) = (\rho, \theta + t)$.

Topology  of $M$ is essentially eliminated from our subsequent  considerations
due to the following fact.

\begin{prop}\label{propweakreduction}
Given any $\epsilon >0$ there exists a $C^{\infty}$ diffeomorphic embedding 
$E: P \to M$ such that  $\mu(E(P))>1-\epsilon$ and  $E_*(\lambda)$ is a scalar multiple of
$\mu$.

\end{prop}
\begin{proof} Let $E_0: P\to M$ be an arbitrary $C^{\infty}$ diffeomorphic 
embedding.
Then $(E_0)_*(\lambda)=\rho\mu$, where $\rho$ is a positive $C^{\infty}$ 
function defined on 
the manifold with boundary $E_0(P)\subset M$ and by definition  
$\int_{E_0(P)}\rho d\mu=1$.
Extend $\rho$ to a positive  $C^{\infty}$ function on the whole 
manifold $M$   which we will still
denote by $\rho$ in such a way that $\int_M\rho d\mu< (1-\epsilon)^{-1}$.
 By the  \cite{Mos} there exists a diffeomorphism $h: M\to M$ such that
 $h_*(\rho \mu)=(\int_M\rho d \mu)\mu$. 
 Now take $E=h\circ E_0$.
 \end{proof}

Thus, it is enough to prove the theorem for the particular case 
of  the open  manifold $P$  with the measure $\lambda$.  A compactly supported diffeomorphism of $P$ is translated to $E(P)\subset M$ via the embedding $E$
and is extended by identity to the rest of $M$.

In fact we will deal with the product $\DD^{m-1}\times S^1$ of  the closed ball with the circle and will construct a diffeomorphism with  the desired minimal set   which on the boundary  $S^{m-2}\times S^1$ coincides with a certain  element  $S_{\alpha}$ of the  action $\mS$. By slightly shrinking the horizontal direction by a homothety
and extending the diffeomorphism to $P$  by  vertical rotations with 
the angle depending smoothly on the radius and decreasing from $\alpha$ 
to identical  zero near the boundary we produce the   desired compactly 
supported diffeomorphism of $P$.

\subsection{Statement and description of the construction.}

\begin{prop} Given any $s \in [0,1)$, there exists a sequence $\a_n$ of
rationals and a sequence $H_n$ of diffeomorphisms preserving $\lambda$ and
constructed as in (\ref{conjugations2}) and (\ref{conjugations3}), such
that the sequence of diffeomorphisms $H_n S_{\a_{n+1}} H_n^{-1}$ converges
in the $C^\infty$ topology to a diffeomorphism of $\DD^{m-1}\times S^1$
preserving $\lambda$ and a compact invariant minimal set of measure $s$. 

Moreover, $f$ can be made arbitrary close to identity in the $C^\infty$
topology. 
 
\end{prop}

\subsubsection{Criterion of minimality}\label{preliminaries} 

\begin{defi} Let $f$ be a map on a complete separable metric 
space $M$. Given $\epsilon > 0$ and a subset $A \subset M$, we say that
$f$ is \emph{ $\epsilon$-minimal on $A$} if given any two points $x, y \in
A$, there exists $n=n(x,y) \in \N$ such that $d(f^n(x),y) < \epsilon$. 

If $m$ is such that $n$  as above can be chosen less than $m$ for any pair
$x, y \in A$ we say that  $f$ is \emph{$(\epsilon, m)$-minimal on $A$.} 
\end{defi}

\begin{lemm}\label{lemmamin} If $A_n \subset M$ is a sequence of closed
sets 
such that
\begin{enumerate}
\item For every $n \in \N$, $f(A_n) = A_n$,
\item There exists a sequence $\epsilon_n \rightarrow 0$ such that $f$ is
$\epsilon_n$-minimal on $A_n$,
\end{enumerate}
 then $\bigcap_{n \in \N} A_n$ is a minimal set for $f$.
\end{lemm}

\begin{proof} We have that $A=\bigcap A_n$ is a closed set and that $f
(A) \subset A$. On the other hand, given any two points $x,y \in A$, there
exists a sequence 
$k_n$ such that $d(f^{k_n}(x),y) < \epsilon_n$, hence $ A$ is a minimal
set for $f$. 
\end{proof}

\begin{rema} In our construction, we will have $A_n$ decreasing and
$\mu(A_n)\rightarrow s$, where $s$ is arbitrarily chosen in $(0,1)$. 
\end{rema}

\subsubsection{ Properties of  the successive conjugating diffeomorphisms
$H_n$.}

 In addition to the  condition (\ref{convergence}) that yields convergence
in (\ref{conjugations}) we will now list certain extra  conditions on the
sequences $\a_n$ and $H_n$ under which the  limiting map will possess a
minimal invariant set of the given measure $s \in [0,1)$.   In the next
subsection we will show how at step $n$ $h_n$ is constructed such that the
required conditions on $H_n$ hold, then choose $\a_{n+1}$  to satisfy the
remaining requirements. 

Denote by $B(r)$   the ball in $\RR^{m-1}$ of radius $r$ around the origin
and let $r_s$ be such that $\lambda ( B ({r_s}) ) =s$. For $n$ large
enough we consider in $M$ the set $C_n = B (r_s + {1 \over n}) \times
S^1$.

\begin{lemm}\label{Lemma1} If $\a_n$ is a sequence of rationals and $H_n$
is a sequence of diffeomorphisms of $M$ satisfying
(\ref{conjugations2}) and (\ref{conjugations3}), and if  for some sequence
$\epsilon_n\to 0$,  one can find  a sequence $m_n$ such that

\begin{enumerate}
\item\label{la)}  $H_n \equiv H_{n-1}$ outside  the interior of
${C}_{n-1}$,
\item\label{lb)} $H_n S_{\a_{n+1}} H_n^{-1}$ is $(\epsilon_n,m_n)$-minimal
on $H_{n} (C_{n})$,
\item\label{lc)} $\displaystyle{ |\a_{n+1} - \a_{n}| <  {1 \over 2^n m_{n-1}
{||H_n||}_{C^n}^n }}$,
\end{enumerate}
 then the sequence $H_n S_{\a_{n+1}}  H_n^{-1}$ converges in the
$C^\infty$ topology to a diffeomorphism $f$ that preserves $\lambda$ and 
is minimal on the set $\bigcap H_n(C_n)$ whose measure is equal to $s$. 
\end{lemm}

\begin{proof} Condition (\ref{lc)}) implies convergence to a volume
preserving $C^{\infty}$ diffeomorphism $f$ (Cf. (\ref{convergence}) in \S
\ref{general scheme}).  We claim that (\ref{la)}) implies that
$f(H_n(C_n)) = H_n(C_n)$: Indeed since $H_{n+1}$ and $H_n$ coincide on the
boundaries of $C_n$ we have that $H_{n+1} (C_n) = H_n(C_n)$ and similarly, because the sequence $C_n$ is nested, $H_l(C_n) = H_n(C_n)$ for any $l \geq n$. Hence $H_l
S_{\a_{l+1}} H_l^{-1} H_n(C_n) = H_{l} S_{\a_{l+1}} C_n = H_l C_n = H_n
(C_n)$ and the claim follows as $l \rightarrow \infty$. The sets
$H_n(C_n)$  are nested,  and $\lambda(H_n(C_n)) = \lambda
(C_n) \rightarrow s$, hence $\lambda ( \bigcap H_n(C_n)
 ) = s$.

Condition (\ref{lc)}) implies as in Lemma \ref{fnn} of \S
\ref{sectionthreemeasures} that $||f^m - H_n \circ S_{ \a_{n+1}}^m \circ
H_n^{-1}|| \leq {m \over 2^n m_n}$, hence (\ref{lb)}) implies that $f$
itself is $(\epsilon_{n} + 1  / 2^n , m_n)$-minimal on $H_{n} (C_{n})$. The
fact that the set $\bigcap H_n(C_n)$ is a minimal set for $f$ then follows
from  Lemma \ref{lemmamin}.
\end{proof}

In the next lemma we will  show how to deduce  conditions (\ref{lb)}) and
(\ref{lc)}) of the previous lemma from certain  conditions on $H_n$ and
$\a_n$ that,  unlike  (\ref{lb)}), do not involve explicit iterations. 

\begin{defi} We will say that a collection $\xi$ of subsets of the  set
$C= B(r) \times S^1$ is {\emph a $\nu$-collection of
horizontal strips} if the  elements $\xi_i$ of $\xi$ are disjoint and  have the  form
$\xi_i = B(r) \times J_i $, where $J_i$ is a  segment of $S^1$  of length
$|J_i| \geq \nu$. 
\end{defi}

\begin{defi}   A set $A \subset B \subset M$ is said to be $\epsilon$-dense in
$B$ if given any $y \in B$ there 
exists $x \in A$ such that $d(x,y) \leq \epsilon$.
\end{defi}

\begin{rema} \label{rrr}
Let $\a_{n+1} = {p_{n+1} /  q_{n+1}}$ with $p_{n+1} $ and $ q_{n+1}$
relatively prime. Then given any $(\rho, \theta) \in C_n$ and any
horizontal strip $\xi$ in $C_n$ of width larger than $1 / q_{n+1}$ there
exists   
 $m \leq q_{n+1}$ such that 
 $$S_{\a_{n+1}}^m(\rho, \theta) = (\rho, \theta + m \a_{n+1} \ {\rm mod }
(1) ) \in \xi.
 $$
\end{rema}

The following lemma implies Lemma \ref{Lemma1} with  $m_n = q_{n+1}$.

\begin{lemm}\label{Lemma2}  If $\a_n = {p_n / q_n} $ is a sequence of
rationals and $H_n$ is a sequence of diffeomorphisms of $M$ satisfying
(\ref{conjugations2}) and (\ref{conjugations3}), and if we have for some
sequence $\epsilon_n\to 0$    a sequence $\xi_n=\{\xi_{n,i}\}$ of ${1 \over q_{n+1}}$-collections of horizontal strips in $C_n$
such that 
\begin{enumerate}
\item \label{zz1} $H_n \equiv H_{n-1}$ outside the interior of ${C}_{n-1}$, 
\item \label{rrb2} The set $\cup_i H_n (\xi_{n,i})$ is $\epsilon_n$-dense in $H_{n} (C_{n})$  
\item \label{rrb} For any element  $\xi_{n,i} \in \xi_{n}$, ${\rm diam} (H_n(\xi_{n,i})) \leq \epsilon_n$, 
\item \label{zz4} $\displaystyle{|\a_{n+1} - \a_{n}| < {1 \over 2^n q_n {\|H_n\|}_{C^n}^n}}$, \end{enumerate}
 then the  conclusions of Lemma \ref{Lemma1} hold.
\end{lemm}

\begin{proof} Fix an arbitrary $z \in H_n(C_n)$. We have from Remark
\ref{rrr} that for any element $\xi_{n,i} \in \xi_n$ there exists $m \leq
q_{n+1}$ such that  $H_n S_{\a_{n+1}}^m H_n^{-1}  (z) \in H_n (\xi_{n,i}).$
It follows then from (\ref{rrb2}) and (\ref{rrb}) that $H_n S_{\a_{n+1}}^m H_n^{-1}$ is $(3 \epsilon_n, q_{n+1})$-minimal on $H_n(C_n)$ which is (\ref{lb)}) of Lemma \ref{Lemma1}. 
Finally,  with our choice of $m_n$  Condition (\ref{lc)}) matches with the corresponding
condition in Lemma \ref{Lemma1}. \end{proof}

\subsection{Construction of $h_n$ and choice of $\a_{n+1}$.} 

 We fix  a sequence $\epsilon_n < {1 \over n^3}$ for all $n \in \NN$.  
Now we assume that the sequence $H_l= h_l \circ ... \circ h_1$ and
$\a_{l+1}$ are constructed up to $l=n-1$ and satisfy the conditions of
Lemma \ref{Lemma2} and proceed to construct $h_n$ and choose $\a_{n+1}$
such that the conditions of the lemma  hold for $l=n$.
 
\subsubsection{Construction of $h_n$.}

Define 
\begin{eqnarray} \epsilon_n':= {\epsilon_n \over
{||H_{n-1}||}_{C^1}}. \label{rrrepsilon} \end{eqnarray} 

Define 

$$\overline{C}_n := 
B \left(r_s + {1 \over n} + {\epsilon_n'}^3 \right) \times S^1.
$$
Clearly, $C_n\subset \overline{C}_n\subset C_{n-1}$,and both inclusions
are strict. Cf. Fig 4.1.

We pick $\nu_n$ sufficiently small and consider a $\nu_n$-collection of horizontal strips $\xi_n$ of $C_n$ such that
\begin{itemize}
\item[(i)] $\xi_n$ is invariant by $S_{1 /q_n}$;
\item[(ii)] For every $\xi_{n,i} \in \xi_n$, $\lambda (\xi_{n,i}) \leq \lambda (\tilde{B} ({\epsilon_n'} / 4))$  (where $\tilde{B} (r)$ denotes the  ball in $\RR^m$ of radius $r$ centered at the origin); 
\item[(iii)] $\lambda(\cup_i(\xi_{n,i})) \geq \lambda (\overline{C}_n) - \lambda(\tilde{B} ({\epsilon_n'} / 4))$.
\end{itemize}


Condition (iii) is possible since $\lambda(\overline{C}_n) \leq \lambda(C_n) + \lambda(\tilde{B} (\epsilon'_n)$. Due to (ii) in the definition of $\xi_n$ it is possible to ask from a volume preserving  diffeomorphism $h_n$ to send each strip $\xi_{n,i}$ inside a ball of radius $\epsilon_n'$ which in turn is satisfied as we choose $\nu_n$ sufficiently small. More precisely we have

\begin{prop} \label{propconjugacy} There exists a diffeomorphism $h_n$ of
class $C^\infty$ on $
\DD^{m-1}\times S^1$ such that 
\begin{enumerate}
\item \label{ee1} $h_n$ is volume preserving,
\item \label{ee2}$h_n \circ S_{\a_n} = S_{\a_n} \circ h_n$,
\item \label{ee3} $h_n \equiv {\rm Id}$ outside $\overline{C}_n$,
\item \label{ee4} For any strip $\xi_{n,i} \in \xi_n$, ${\rm diam} (h_n
(\xi_{n,i})) \leq \epsilon_n'$.
\end{enumerate}
\end{prop}

Condition (\ref{ee2}) can be satisfied  since the collection $\xi_n$ is
invariant by $S_{\a_n}$ and (\ref{ee1}) and  (\ref{ee4}) are compatible
since we have imposed that $\lambda (\xi_{n,i}) \leq  \lambda (B
({\epsilon_n'} / 2))$. 

\vspace{1cm}

 \setlength{\unitlength}{0.00050000in}
\begingroup\makeatletter\ifx\SetFigFont\undefined%
\gdef\SetFigFont#1#2#3#4#5{%
  \reset@font\fontsize{#1}{#2pt}%
  \fontfamily{#3}\fontseries{#4}\fontshape{#5}%
  \selectfont}%
\fi\endgroup%
{\renewcommand{\dashlinestretch}{30}
\begin{picture}(7749,7002)(0,-10)
\path(12,6975)(12,1725)
\path(3012,6975)(3012,1725)
\path(312,6975)(312,1725)
\path(2712,6975)(2712,1725)
\path(612,6975)(612,1725)
\path(2412,6975)(2412,1725)
\path(612,2100)(2412,2100)
\path(612,2175)(2412,2175)
\path(612,2475)(2412,2475)
\path(612,1800)(2412,1800)
\path(612,2550)(2412,2550)
\path(612,2850)(2412,2850)
\path(612,1425)(2412,1425)
\path(312,900)(2712,900)
\path(87,375)(3087,375)
\path(612,1500)(612,1425)
\path(2412,1500)(2412,1425)
\path(312,975)(312,900)
\path(2712,975)(2712,900)
\path(87,450)(87,375)
\path(3087,450)(3087,375)
\thicklines
\path(5037,3000)(5037,3225)
\path(5037,2100)(5037,3225)(5787,3225)
	(5787,2100)(5937,2100)(5937,3225)
	(6687,3225)(6687,2100)(6837,2100)
	(6837,3225)(7512,3225)(7512,2100)
	(7662,2100)(7662,3300)(4962,3300)
\thinlines
\path(7587,1875)(5637,1875)(5637,3075)
	(5712,3075)(5712,1950)(6537,1950)
	(6537,3000)(6612,3000)(6612,1950)
\path(6612,1950)(7362,1950)(7362,3000)
	(7437,3000)(7437,1950)(7737,1950)
	(7737,3375)(5637,3375)
\path(7587,3375)(5637,3375)(5637,4575)
	(5712,4575)(5712,3450)(6537,3450)
	(6537,4500)(6612,4500)(6612,3450)
\thicklines
\path(5037,2100)(5037,3225)(5787,3225)
	(5787,2100)(5937,2100)(5937,3225)
	(6687,3225)(6687,2100)(6837,2100)
	(6837,3225)(7512,3225)(7512,2100)
	(7662,2100)(7662,3300)(4962,3300)
	(4962,4725)(5787,4725)(5787,3600)
	(5937,3600)(5937,4725)(6687,4725)
	(6687,3600)(6837,3600)(6837,4725)
	(7512,4725)(7512,3600)(7662,3600)
	(7662,4800)(4962,4800)
\thinlines
\path(6612,3450)(7362,3450)(7362,4500)
	(7437,4500)(7437,3450)(7737,3450)
	(7737,4875)(6912,4875)(6912,4875)
\put(1437,3000){\makebox(0,0)[lb]{\smash{{{\SetFigFont{8}{9.6}{\familydefault}{\mddefault}{\updefault}.}}}}}
\put(1437,3150){\makebox(0,0)[lb]{\smash{{{\SetFigFont{8}{9.6}{\familydefault}{\mddefault}{\updefault}.}}}}}
\put(1137,2250){\makebox(0,0)[lb]{\smash{{{\SetFigFont{8}{9.6}{\familydefault}{\mddefault}{\updefault}$\xi_2$}}}}}
\put(1137,1875){\makebox(0,0)[lb]{\smash{{{\SetFigFont{8}{9.6}{\familydefault}{\mddefault}{\updefault}$\xi_1$}}}}}
\put(1137,2625){\makebox(0,0)[lb]{\smash{{{\SetFigFont{8}{9.6}{\familydefault}{\mddefault}{\updefault}$\xi_3$}}}}}
\put(1437,3300){\makebox(0,0)[lb]{\smash{{{\SetFigFont{8}{9.6}{\familydefault}{\mddefault}{\updefault}.}}}}}
\put(1137,1050){\makebox(0,0)[lb]{\smash{{{\SetFigFont{10}{12.0}{\familydefault}{\mddefault}{\updefault}$C_n$}}}}}
\put(1062,525){\makebox(0,0)[lb]{\smash{{{\SetFigFont{10}{12.0}{\familydefault}{\mddefault}{\updefault}$\overline{C}_n$}}}}}
\put(1062,0){\makebox(0,0)[lb]{\smash{{{\SetFigFont{10}{12.0}{\familydefault}{\mddefault}{\updefault}$C_{n-1}$}}}}}
\put(5862,1050){\makebox(0,0)[lb]{\smash{{{\SetFigFont{10}{12.0}{\familydefault}{\mddefault}{\updefault}$h_n(C_n)$}}}}}
\put(7437,2550){\makebox(0,0)[lb]{\smash{{{\SetFigFont{8}{9.6}{\familydefault}{\mddefault}{\updefault} }}}}}
\put(7437,2475){\makebox(0,0)[lb]{\smash{{{\SetFigFont{8}{9.6}{\familydefault}{\mddefault}{\updefault}         ...}}}}}
\put(5037,3900){\makebox(0,0)[lb]{\smash{{{\SetFigFont{7}{8.4}{\familydefault}{\mddefault}{\updefault}$h_n(\xi_i)$}}}}}
\end{picture}
}

$$ \hspace{3.8cm} {\rm Fig. \ 4.1} \hspace{4.2cm} $$

\vspace{1cm}

\subsubsection{ Checking the conditions of Lemma \ref{Lemma2} and choice
of $\a_{n+1}$.}   By  ({\ref{ee1}) and (\ref{ee2}) we can define the
conjugating diffeomorphisms $H_n = h_1 \circ ... \circ h_n$ as in
({\ref{conjugations2}) and (\ref{conjugations3}). By definition of
$\overline{C}_n$,  (\ref{ee3}) implies (\ref{zz1}) of Lemma \ref{Lemma2}. 
Due to our choice of $\epsilon'_n$ in (\ref{rrrepsilon}) 
Condition ({\ref{ee4}) implies (\ref{rrb}) of Lemma \ref{Lemma2}.
Next, due to (iii) in the definition of $\xi_n$ we have that any ball of radius $\epsilon_n'$ centered at a point in $\overline{C}_n$ intersects some $h_n(\xi_{n,i}) \subset \overline{C}_n$. Hence $\cup_i h_n (\xi_{n,i})$ is $\epsilon_n'$-dense in $\overline{C}_n = h_n (\overline{C}_n)$. Again, from our choice of $\epsilon_n'$ (\ref{rrb2}) of Lemma \ref{Lemma2} follows.

Finally, we choose  
$\a_{n+1} = {p_{n+1} / q_{n+1}}$ such that (\ref{zz4}) in Lemma
\ref{Lemma2} is satisfied and 
such that  $q_{n+1} \geq {1 / \nu_n}$ in order for the $\nu_n$-collection $\xi_n$ considered  to be a $
{1 /  q_{n+1}}$-collection as required in Lemma \ref{Lemma2}.

All the conditions of Lemma \ref{Lemma2} being satisfied by the sequences
$\a_n$ and $H_n$ thus constructed, the diffeomorphisms  $H_n S_{\a_{n+1}}
H_n^{-1}$ converge in the $C^\infty$ topology to a diffeomorphism with  a
minimal set  of measure $s$. \carre

\subsubsection{Connectedness and rigidity}
The set $C$  is connected since it is  the intersection of the nested sequence of connected sets
$H_n(C_n)$. Furthermore, in the two--dimensional case $m=2$ the complement to each set $H_n(C_n)$ to  $\DD \times S^1$ consists of two  open connected components and  this remains true  for $C$.   Each connected component of the complement is dense in $C$.
For  $m>2$ the set $\DD^{m-1}\times S^1\setminus C$ is connected.

 It follows immediately from our construction that the  diffeomorphism $f$ is  $C^{\infty}$ \emph{rigid}, 
namely $f^{q_n} \to\Id$  as $n\to\infty$ in  the $C^{\infty}$ topology.

\subsubsection{Special choice of $h_n$ and $\a_{n+1}$}\label{sbsspecial}
As was mentioned in Section \ref{almostallall} the most natural way to ensure commutativity relation (\ref{lb)}) in Proposition \ref{propconjugacy} is to take 
the fundamental  domain $\Delta_n = \DD^{m-1}\times[0,  {1\over{q_n}}]$,  construct a volume preserving diffeomorphism $g$ of $\Delta$ onto itself, identical near its boundary and satisfying
conditions (\ref{ee3}) and \ref{ee4}) and  extend $g$ by periodicity. It is also  natural to choose $q_{n+1}$ to be  a multiple of $q_n$.

For the purposes of the further discussion we will call  the choice of  $h_n$'s and  $q_n$'s  described above simply {\em special}.

With these  special choices one produces  a nested sequence of partitions invariant under
 $H_n S_{\a_{n+1}} H_n^{-1}$.
Those partitions are $H_n(\zeta_n)$ where $\zeta_n$ is simply the partition into the ``layers''
$ \DD^{m-1}\times[{k\over{q_n}},  {k+1\over{q_n}}], \,\,k=0,\dots,q_n-1.$ This layers are simply images of the fundamental domain $\Delta_n$ under iterates of the rotation $S_{\a_n}$.
Restricted to the minimal set $C$ this sequence of  partitions  becomes exhaustive in the measurable sense (converges to $\epsilon$ in the more customary terminology of ergodic theory) due to Lemma \ref{Lemma2} and \ref{propconjugacy}, (\ref{ee4}).
Notice also that  on the boundaries of the partition elements  the conjugacies $H_n$ stabilize:
For the boundary of $\zeta_n$ all $h_m$ with  $m>n$ are identities.
This implies in particular that the set $C$  contains  a  dense set of smooth  $m-1$--dimensional discs 
\begin{eqnarray}\label{specialdiscs}
D_{n,k}=:H_n(B(r_s)\times\{{k\over{q_n}}\}),\,\, n=1,2,\dots; k=0,1,\dots q_n-1.
\end{eqnarray}

The limit partition $\xi=:\lim_{n\to\infty}H_n(\zeta_n)$ is invariant under $f$.
Notice that  elements of that partition are closed connected sets.
Furthermore,  by our construction not only $f$ but every power of it  is minimal.

\begin {lemm} Every element  $c\in\xi$ is nowhere dense in $C$  and $f^n c\neq c$.
\end{lemm}

\begin{proof}
Suppose $f^n c= c$. Then the boundary of $c$ is nowhere dense and  the union of its images is
a nowhere dense $f$ invariant sucset of $C$ contradicting minimality.
If $c$ has an interior  then all  images  of it are disjoint and their union is an invariant open set.
By minimality it coincides with  $C$ and hence  $c$ is both open and closed, contradicting 
connectedness of $C$

\end{proof}

\begin{lemm}
The discs $D_{n,k}$ (and hence their images) are elements of the partition $\xi$.
\end{lemm}

\begin{proof} These discs are parts of larger discs $H_n(\DD^{m-1}\times S^1)$; hence 
 the  discs $D_{n,k}$  are ordered cyclically, in other words, the same way as the numbers
 ${k\over{q_n}},\,\,(\mod 1)$. Furthermore, one can consider more discs of the form $f^{q_m}D_{n,k},\,\, m\in \ZZ$ and all those discs are ordered the same  way as 
 $n,k$, such that ${{k\over{q_n}}\le  y \le {{k+1}\over{q_n}}}$. Obviously
 every disc  $D_{n,k}$  lies within an element of $\xi$  which we will denote $C_{n,k}$ and every image $f^ lD_{n,k}$  for $l\neq 0$ belongs to a different  element. But since both $f^{q_n}$ and $f^{-q_n}$ converge to identity as $n\to\infty$ one concludes that  $C_{n,k}=D_{n,k}$ because the whole element $C_{n,k}$ must lie between the discs $f^{q_n}D_{n,k}$ and $f^{-q_n}D_{n,k}$.
  \end{proof}

To finish the proof of  Theorem \ref{theoremmain} it is sufficient to prove the following statement.
 
 \begin{prop} With any  special choice of  diffeomorphisms $h_n$ and $\a_n$  
  for any $x\in C$ and any sufficiently small open neighborhood $V$ of $x$, the
intersection $C\cap V$ has topological dimension $m-1$ and is not homeomorphic  to the 
product of a disc with  a Cantor set. 
 \end{prop}
 
 \begin{proof}
 The union of  $m-1$-dimensional discs  $D_{n,k}\subset C$ (\ref{specialdiscs}) is dense in $C$.
  Now assume that  for some $x\in C$ there is an open neighborhood
 $V$ such that $C\cap V$ is homeomorphic to the product of $\RR^{m-1}$ with a Cantor set $K$.
 Since elements of $\xi$ are connected  for any $x\in C\cap V$  the  fiber $F_x$ in the product containing $x$  belongs to  a single element of $\xi$. By rigidity $f^{q_n}(F_x)$ converge to $F_x$ and hence lie on different fibers. Now take the  fiber $F_x$ corresponding  to an endpoint  of the complimentary interval  to the Cantor set.  Replacing $V$ by a slightly smaller compact set which we still denote by $V$ we thus obtain the  image  $V$
 of $[-1,1]\times[0,1]$  under a homeomorphism  which we will denote $\mathcal H$ such that $F_x=\mathcal H (\{0\}\times[0,1])$ and  $\mathcal H ((0,1]\times[0,1])$  does not intersect $C$ while there is a sequence of negative numbers
 $x_n\to 0$ such that   $\mathcal H(F_{f^{q_n}x})\subset \{x_n\}\times[0,1]$. Since $f^{q_n}\to \Id$
 we may assume without loss of generality that  $\mathcal H(\{x_n\}\times[0,1])\subset  C$.
Now consider the sets $f^{q_n}\mathcal H ((0,1]\times[\epsilon,1-\epsilon])$
By invariance of $C$  these sets are disjoint from $C$.
Since $f^{q_n}\to \Id$ for a sufficintely small $\epsilon>0$ and any large enough $n$ this  set 
lies in $V$ and hence $\mathcal H^{-1}f^{q_n}\mathcal H ((0,1]\times[\epsilon,1-\epsilon])\cap
(\{x_m\}\times[0,1])=\emptyset $ for any sufficintly large $m$.  Consequently $\mathcal H^{-1}f^{q_n}\mathcal H ([0,1]\times[\epsilon,1-\epsilon])\subset [0,1]\times [0,1]$. But this implies that
$F_{f^{q_n}}f\subset F_x$, a contradiction.
\end{proof}

\begin{prop}\label{propNakayama} In dimension two for the special  choice of parameters the  set $C$ is never locally connected.
\end{prop}

\begin{proof}  Every  complete locally connected space is locally path connected. As is shown in \cite{BNW}  such a set must be  a simple closed
curve.  This is impossible since for exmaple,  the  closed arcs   $f^{q_n}D_{1,0}$ are disjoint  due to the ordering and
due to rigidity their  diameters are bounded away from zero. 
\end{proof}

\subsection{Other versions of the construction}Notice that  in Handel's construction \cite{Ha}  for the two--dimensional case  the minimal set is a pseudo-circle.
Since a pseudo-circle does not contain any embedded  topological arcs  any set $C$ 
obtained with a special choice of parameters as in Section \ref{sbsspecial} is  different from 
those obtained by Handel \cite{Ha} which are  pseudocircles.

 Handel's construction may be put into our framework to  make  the
set of positive (and, in fact, given) measure. The main  mechanism that 
makes the resulting set nontrivial  in that construction is a  twist of the standard fundamental   domains along the orbits of the rotation  rather than 
spreading the orbits  within those  domains as for    special choices  of parameters. 

For another special choice of parameters one can actually make $C$ a simple closed curve of positive  area.  In this case  the diffeomorphisms $h_n$ converge  to identity in $C^0$ topology
and  their products $H_n$  converge to  an  area--preserving continuous map which is not injective.

\begin{prob}
Let   $C$ be a minimal set  for a diffeomorphism  $f$
of  a   two--dimensional  manifold $M$.
\begin{enumerate}
\item Can $C$ be nowhere dense and path--connected?
\item Can $C$ be locally homeomorphic to the product of the line and the Cantor set?
\end{enumerate}

\end{prob}

\subsection{Ergodic properties of the  minimal set}
The restriction of the diffeomorphism $f$ to  the minimal set $C$ preserves the Lebesgue measure,
or, more precisely, the absolutely continuous  probability measure   $\lambda_C$ obtained by  normalizing the restriction of the measure $\lambda$ to $C$. Our construction guarantees unique ergodicity of  $f$ restricted to the set.
$C$. Furthermore, the special choice of parameters described in  Section \ref{sbsspecial}
makes $(f,\lambda_C)$ metrically isomorphic to the rotation of the circle by the angle 
$\a=\lim_{n\to{\infty}}\a_n$ as in the construction of \cite[Section 4]{AK}. On the other hand,  one can choose parameters in a different way to make $f$ weakly mixing with respect to $\lambda_C$. Other modifications are also possible, see 
\cite{AK}. 

On the other hand, by modifying our construction along the lines described in \cite{W} 
 one may violate unique ergodicity;
for example make $\lambda_C$ decompose into any given finite or countable number of ergodic components. Furthermore, one may modify our construction by making  the  transformations $h_n$
not volume preserving in such a way that  the minimal set will still have   positive (and prescribed)
Lebesgue measure but invariant measure(s) of $f$ supported by the set $C$ will be all singular or some of them will be absolutely continuous and some singular.
Presenting all these variations in detail will take too much space and will look rather tedious.

\begin{prob}
Is it possible for a minimal set of  a diffeomorphism of a two--dimensional manifold to carry 
an invariant mixing measure?
\end{prob}

Notice that this  question is open both for nowhere dense invariant  minimal sets and for minimal diffeomorphisms of $\TT^2$.

\section{Transitive non ergodic diffeomorphisms.}\label{sectiontrnonerg}

Let $f$ be a map on a complete separable metric space $M$. Denote 
by ${\phi}_f \subset M$ the set of points with dense orbit, i.e. $x 
\in {\phi}_f$  if and only if $\overline{ \lbrace f^n(x) \rbrace_{n 
\in \Z}} = M$. Recall that a map $f$  is said to be topologically transitive 
(or simply transitive)  if ${\phi}_f$ is not empty and in this case it is easy to see 
that ${\phi}_f$ is a $G_\delta$ dense set in $M$. When $M$ is a 
Riemannian manifold it is natural to ask how small can the measure of 
${\phi}_f$ be for a topologically  transitive map $f$, especially in the case of a 
volume preserving map $f$.  If such a map  is ergodic then the set  $\phi_f$ has full measure. The same may be true for non-ergodic  maps, see e.g.  \cite{furstenberg}. We show in this section that the opposite  can also be true. To build examples,
 we again use  an appropriately adjusted version of the approximation  by conjugation method  of  Section \ref{general scheme}.

\begin{theo}\label{yytheorem} Let $M$ be an $m$-dimensional 
differentiable manifold with a nontrivial circle action
  ${\mS} = {\lbrace S_t \rbrace }_{t \in\RR}, \  S_{t+1}=S_t$ preserving a smooth volume 
$\mu$. There exists a $C^{\infty}$ diffemorphism $f$ of $M$ 
preserving $\mu$ such that $f$ is transitive and $\mu ( {\phi}_f) =0$.
Furthermore, $f$ can be chosen arbitrary close to identity in the 
$C^\infty$ topology.
\end{theo}

\subsection{Reduction to the case 
of ${\DD}^{m-1} \times S^1 $}
We need a  statement  somewhat similar to but stronger than  Proposition \ref{propweakreduction}. Its proof  goes along the same lines as proofs of similar statements  in \cite{AK, W}. 
Let $\lambda$ be the product of Lebesgue measures on ${\DD}^{m-1}$ and $S^1$. Denote by
$\mathcal R$ the standard ``vertical'' action of  $S^1$ on  ${\DD}^{m-1}\times S^1$.

\begin{prop}\label{reductiondirect}Let $M$ be an $m$-dimensional 
differentiable manifold with an effective circle action
  ${\mS} = {\lbrace S_t \rbrace }_{t \in\RR}, \  S_{t+1}=S_t$ preserving a smooth volume 
$\mu$. Let $B=:\partial M\cup F\cup(\underset{q}\bigcup F_q)$ (see section \ref{general scheme}).  There exists a continuous surjective map $F: {\DD}^{m-1} \times S^1 \to M$  with the following properties
\begin{enumerate}
\item The restriction of $F$ to the interior $ {\BB}^{m-1} \times S^1$ is a $C^{\infty}$ diffeomorphic embedding;
\item $\mu(F(\partial( {\DD}^{m-1} \times S^1))=0$;
\item $F(\partial( {\DD}^{m-1} \times S^1)\supset B$;
\item$F_*(\lambda)=\mu$;
\item $\mS\circ F=F\circ\mathcal R$.
\end{enumerate}
\end{prop}

Now suppose $f: {\DD}^{m-1} \times S^1 \to  {\DD}^{m-1} \times S^1 $ is a topologically transitive diffeomorphism  such that   $\mu ( {\phi}_f) =0$, suppose furthermore that $f | \partial( {\DD}^{m-1} \times S^1=\mathbb R_{\alpha}$ for some $\alpha$. Then one can define  a homeomorphism
$g: M\to M$  as $$g(x)=F^{-1}(f(Fx))\,\, \text{for}\,\,  x\in F({\BB}^{m-1} \times S^1),\,\,\text{and}$$ 
$$g(x)= \mathbb R_{\alpha}(x)\,\, \text{for}\,\,\in F(\partial( {\DD}^{m-1}) \times S^1).$$. If moreover,   $f$ has the same jets of all orders on  as $\mathbb R_{\alpha}(x)$ on  $\partial( {\DD}^{m-1}) \times S^1)$ and the difference is sufficiently flat (i.e. the jets of the difference decay fast enough to near the boundary; this condition depends on $F$) then $g$ constructed  as above is a $C^{\infty}$ diffeomorphism  satisfying the assertions of Theorem \ref{yytheorem}.

\subsection{Construction of the diffeomorphism $h_n$.} We will describe the construction only for the case $m=2$. This makes it easier to visualize (see Fig 5.1. below). In the general case  the two
``bumps'' have to be replaces by a curve whose projection to the horizontal  $\DD^{m-1} $  direction is  almost dense. For example, for $m=3$ one may take a spiral beginning at the center of the disc, going outside almost to the boundary, then changing direction and  returning to the center.  Closeness to the rotation with  all derivatives   near the boundary required for the reduction described above can be guaranteed  by taking the  successive conjugating diffeomorphisms sufficiently flat near the boundary. 

Assume that the diffeomorphisms $h_l$ have been constructed 
up to $l=n-1$ and that the rational $\a_n = p_n / q_n$ has been 
chosen.

\begin{defi} \label{yydefinition} For every integer $n \geq 1$ and 
every integer $k \leq q_n - 1$ we define

\begin{eqnarray*}
B_{n,k}^r &=& \left[ {1 \over 2} + {k \over 2 q_n}, {1 \over 2} + {k 
+ 1 \over 2 q_n} \right] \times \left[ 0, {1 \over q_n} \right], \\
  B_{n,k}^l &=& \left[{1 \over 2} - {k + 1 \over 2 q_n}, {1 \over 2} - 
{k \over 2 q_n} \right] \times \left[ 0, {1\over q_n} \right], \\
  I_{n,k}^r &=& \lbrace {1 \over 2} \rbrace \times  \left[{k \over 
q_n^3} + {1 \over 4 q_n^3}, {k \over q_n^3} + {3 \over 4 q_n^3} 
\right] \subset \lbrace {1 \over 2} \rbrace \times [0, {1 \over q_n^2}], \\
I_{n,k}^l &=& \lbrace {1 \over 2} \rbrace \times  \left[ { 1 \over 
q_n^2} + {k \over q_n^3} + {1 \over 4 q_n^3}, {1 \over q_n^2} + {k 
\over q_n^3} + {3 \over 4 q_n^3}  \right] \subset \lbrace {1 \over 2} \rbrace \times [{1 \over q_n^2}, {2 \over q_n^2}]. \\
\end{eqnarray*}

For every integer $n \geq 1$ and for any $\rho_0 \in [0,1]$ define

\begin{eqnarray*}
D_{\rho_0}^r &=& \lbrace (\rho, \theta) \in M \ / \ \rho_0 \leq \rho 
\leq 1, \ \theta \in \TT \rbrace, \\
D_{\rho_0}^l &=& \lbrace (\rho, \theta) \in M \ / \ 0 \leq \rho \leq 
\rho_0, \ \theta \in \TT \rbrace. \\
\end{eqnarray*}

Finally for any $n \geq 1$ define
\begin{eqnarray*}
M_n = \left[ 0,1 \right] \times \left[ {1 \over q_n^{3 \over 2}}, {1 \over q_n} \right].
\end{eqnarray*}

\end{defi}

\begin{prop}[Properties of $h_n$] \label{yyhn} For every integer $n 
\geq 1$, there exists a diffeomorphism $h_n$ with the following 
properties

\begin{enumerate}
\item \label{yyh1} $h_n$ preserves $\mu$,

\item \label{yyh2} $h_n S_{1 \over q_n} = S_{1 \over q_n} h_n$,

\item \label{yyh3} $h_n \equiv {\rm Id}$ on the set $M_n$, in 
particular $h_n(1 / 2 , 0) = (1/2, 0)$,

\item \label{yyh4} For any $k \leq q_n -1$, $h_n(I_{n,k}^j) \subset 
B_{n,k}^j$, $j= r, l$,

\item \label{yyh5} For any $\rho_0 > 1 / n$ we have

\begin{eqnarray*}
h_n \left(D^r_{{1 \over 2} + \rho_0} \right) &\subset&
  D^r_{{1 \over 2} + \rho_0 - {1 \over 2^n}}, \\
h_n  \left(D^l_{{1 \over 2} - \rho_0} \right) &\subset&
  D^l_{ {1 \over 2} - \rho_0 + {1 \over 2^n}}.
\end{eqnarray*}

\end{enumerate}

\end{prop}

\vspace{0.2cm}

\setlength{\unitlength}{0.00066667in}
\begingroup\makeatletter\ifx\SetFigFont\undefined%
\gdef\SetFigFont#1#2#3#4#5{%
  \reset@font\fontsize{#1}{#2pt}%
  \fontfamily{#3}\fontseries{#4}\fontshape{#5}%
  \selectfont}%
\fi\endgroup%
{\renewcommand{\dashlinestretch}{30}
\begin{picture}(5937,2439)(0,-10)
\path(3750,912)(3750,12)
\path(1425,912)(1425,12)
\path(1050,2037)(1050,12)
\path(1350,912)(1500,912)
\path(1350,12)(1500,12)
\path(975,12)(1125,12)
\path(4425,87)(4425,12)
\path(3750,912)(3750,2037)
\path(5925,2037)(5925,2412)
\dashline{60.000}(1650,912)(5925,912)
\path(5925,2037)(4500,912)
\path(5925,1662)(4575,612)
\path(4500,912)(4425,837)
\path(4950,537)(5925,1287)
\path(4425,87)(4500,162)
\path(5925,912)(5400,537)
\path(5925,612)(5700,462)
\path(5100,12)(5550,312)
\path(5925,237)(5625,12)
\path(5025,237)(4725,12)
\path(1875,2337)(5925,2337)
\path(1875,2337)(1800,2337)
\path(1725,2412)(1725,12)(5925,12)(5925,2037)
\path(4425,2037)(4425,2337)
\path(5850,2337)(4425,1287)
\path(5400,2337)(4425,1662)
\path(4950,2337)(4425,2037)
\path(1800,2337)(1725,2337)
\path(975,2337)(1125,2337)
\path(1050,2337)(1050,2037)
\path(3750,2037)(3750,2337)
\path(3750,537)(3750,536)(3750,533)
	(3751,526)(3752,517)(3754,506)
	(3757,497)(3760,489)(3764,483)
	(3769,478)(3775,474)(3783,471)
	(3792,469)(3805,467)(3820,465)
	(3839,464)(3859,463)(3879,463)
	(3892,462)(3899,462)(3900,462)
\path(3750,87)(3750,89)(3751,101)
	(3752,118)(3754,132)(3758,143)
	(3763,149)(3769,154)(3780,158)
	(3794,160)(3811,161)(3823,162)(3825,162)
\path(3660,837)(3658,837)(3653,837)
	(3643,837)(3628,837)(3606,837)
	(3578,837)(3542,836)(3498,836)
	(3447,836)(3389,836)(3325,835)
	(3255,835)(3180,834)(3101,834)
	(3020,833)(2937,832)(2854,831)
	(2771,830)(2690,829)(2610,828)
	(2533,827)(2459,826)(2389,825)
	(2322,824)(2260,822)(2202,821)
	(2148,819)(2098,818)(2053,816)
	(2012,814)(1975,812)(1942,811)
	(1913,808)(1888,806)(1866,804)
	(1848,801)(1834,799)(1823,796)
	(1816,793)(1811,790)(1810,787)
	(1811,784)(1816,780)(1823,776)
	(1834,773)(1848,768)(1866,764)
	(1888,759)(1913,754)(1942,749)
	(1975,744)(2012,738)(2053,732)
	(2098,725)(2148,718)(2202,711)
	(2260,703)(2322,695)(2389,686)
	(2459,678)(2533,669)(2610,659)
	(2690,650)(2771,640)(2854,630)
	(2937,620)(3020,610)(3101,601)
	(3180,592)(3255,583)(3325,575)
	(3389,568)(3447,561)(3498,555)
	(3542,550)(3578,546)(3606,543)
	(3628,541)(3643,539)(3653,538)
	(3658,537)(3660,537)
\path(3750,912)(3750,911)(3750,908)
	(3749,901)(3748,892)(3746,881)
	(3743,872)(3740,864)(3736,858)
	(3731,853)(3725,849)(3717,846)
	(3708,844)(3695,842)(3680,840)
	(3661,839)(3641,838)(3621,838)
	(3608,837)(3601,837)(3600,837)
\path(3750,462)(3750,464)(3749,476)
	(3748,493)(3746,507)(3742,518)
	(3738,524)(3733,528)(3726,531)
	(3717,533)(3706,535)(3693,536)
	(3682,537)(3676,537)(3675,537)
\path(4425,87)(4424,87)(4424,88)
	(4422,88)(4421,89)(4418,89)
	(4415,90)(4412,92)(4408,93)
	(4404,95)(4400,97)(4396,99)
	(4391,101)(4387,103)(4383,106)
	(4379,108)(4375,111)(4372,114)
	(4370,117)(4368,119)(4367,123)
	(4367,126)(4368,129)(4370,133)
	(4373,136)(4377,140)(4383,144)
	(4391,148)(4401,153)(4412,157)
	(4425,162)(4449,169)(4477,177)
	(4505,183)(4535,189)(4566,195)
	(4596,200)(4627,205)(4658,210)
	(4690,214)(4720,218)(4750,222)
	(4778,226)(4804,229)(4826,231)
	(4844,234)(4858,235)(4867,236)
	(4873,237)(4875,237)
\path(4875,237)(5550,312)
\path(3829,462)(3831,462)(3836,462)
	(3846,462)(3861,462)(3883,462)
	(3911,462)(3947,461)(3991,461)
	(4042,461)(4100,461)(4164,460)
	(4234,460)(4309,459)(4388,459)
	(4469,458)(4552,457)(4635,456)
	(4718,455)(4799,454)(4879,453)
	(4956,452)(5030,451)(5100,450)
	(5167,449)(5229,447)(5287,446)
	(5341,444)(5391,443)(5436,441)
	(5477,439)(5514,437)(5547,436)
	(5576,433)(5601,431)(5623,429)
	(5641,426)(5655,424)(5666,421)
	(5673,418)(5678,415)(5679,412)
	(5678,409)(5673,405)(5666,401)
	(5655,398)(5641,393)(5623,389)
	(5601,384)(5576,379)(5547,374)
	(5514,369)(5477,363)(5436,357)
	(5391,350)(5341,343)(5287,336)
	(5229,328)(5167,320)(5100,311)
	(5030,303)(4956,294)(4879,284)
	(4799,275)(4718,265)(4635,255)
	(4552,245)(4469,235)(4388,226)
	(4309,217)(4234,208)(4164,200)
	(4100,193)(4042,186)(3991,180)
	(3947,175)(3911,171)(3883,168)
	(3861,166)(3846,164)(3836,163)
	(3831,162)(3829,162)
\drawline(5550,312)(5550,312)
\path(5325,537)(5326,537)(5330,537)
	(5340,537)(5356,536)(5376,536)
	(5397,535)(5418,534)(5438,533)
	(5456,531)(5475,529)(5493,527)
	(5513,524)(5528,522)(5544,520)
	(5562,517)(5580,514)(5599,510)
	(5619,506)(5639,502)(5659,498)
	(5679,493)(5698,489)(5716,484)
	(5733,480)(5749,475)(5763,471)
	(5776,466)(5788,462)(5793,460)
	(5799,457)(5804,455)(5809,452)
	(5813,449)(5817,447)(5821,444)
	(5824,441)(5827,439)(5829,436)
	(5831,433)(5833,430)(5834,427)
	(5834,424)(5835,422)(5834,419)
	(5834,416)(5832,413)(5831,410)
	(5829,408)(5826,405)(5824,402)
	(5820,400)(5817,397)(5813,394)
	(5809,392)(5805,389)(5800,387)
	(5791,383)(5780,378)(5768,374)
	(5754,369)(5737,364)(5718,358)
	(5696,351)(5672,345)(5647,338)
	(5621,331)(5597,324)(5577,319)
	(5563,315)(5554,313)(5551,312)(5550,312)
\path(5325,537)(5322,537)(5315,537)
	(5303,538)(5284,538)(5258,539)
	(5226,540)(5188,542)(5145,543)
	(5099,546)(5052,548)(5004,551)
	(4956,554)(4911,558)(4867,562)
	(4827,566)(4789,571)(4754,576)
	(4723,581)(4693,587)(4667,594)
	(4642,601)(4620,609)(4599,617)
	(4580,627)(4563,637)(4543,650)
	(4525,664)(4508,680)(4492,696)
	(4477,714)(4463,733)(4450,753)
	(4439,775)(4428,797)(4418,819)
	(4410,842)(4403,865)(4396,889)
	(4391,912)(4386,935)(4383,957)
	(4380,979)(4378,999)(4377,1019)
	(4376,1039)(4375,1057)(4375,1074)
	(4375,1095)(4376,1115)(4377,1134)
	(4378,1154)(4380,1173)(4382,1192)
	(4384,1212)(4387,1231)(4389,1251)
	(4392,1271)(4395,1291)(4398,1311)
	(4401,1331)(4403,1352)(4406,1372)
	(4408,1393)(4410,1415)(4413,1437)
	(4414,1455)(4415,1473)(4416,1493)
	(4417,1514)(4418,1536)(4419,1561)
	(4420,1587)(4421,1616)(4421,1648)
	(4422,1682)(4422,1718)(4423,1756)
	(4423,1795)(4424,1835)(4424,1874)
	(4424,1912)(4425,1945)(4425,1975)
	(4425,1998)(4425,2016)(4425,2028)
	(4425,2034)(4425,2037)
\put(1125,387){\makebox(0,0)[lb]{\smash{{{\SetFigFont{10}{12.0}{\familydefault}{\mddefault}{\updefault}${1 \over q_n^2}$}}}}}
\put(0,462){\makebox(0,0)[lb]{\smash{{{\SetFigFont{10}{12.0}{\familydefault}{\mddefault}{\updefault}.}}}}}
\put(525,1287){\makebox(0,0)[lb]{\smash{{{\SetFigFont{10}{12.0}{\familydefault}{\mddefault}{\updefault}${1 \over q_n}$}}}}}
\end{picture}
}

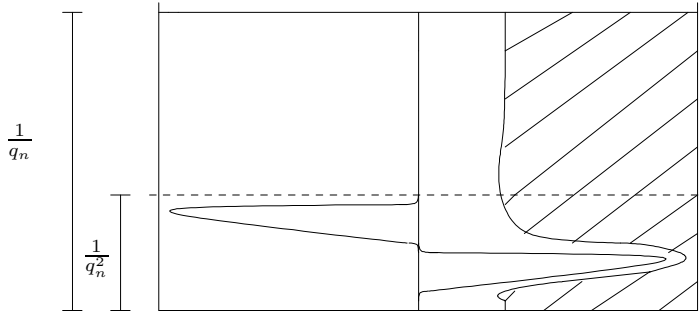
\begin{center} Fig 5.1. \   $h_n \left( \lbrace {1\over 2} \rbrace \times [0, {1 
\over q_n}] \right)$ and $h_n(D^r_{{1 \over 2} + \rho_0})$.

\end{center}

\vspace{0.2cm}

\subsection{Choice of $\a_{n+1}$ and convergence.} \label{yyconvergence}

As in the other sections we choose $\a_{n+1}$ to insure convergence 
of $f_n = H_n S_{\a_{n+1}} H_n^{-1}$ to a diffeomorphism $f$ and have 
in addition that  for any $m \leq q_{n+1}$
\begin{eqnarray} \label{yy1} ||f_n^m - f^m|| \leq {1 \over  n}. \end{eqnarray}

\subsection{Transitivity.} In this paragraph we will prove that the 
limit diffeomorphism $f$ is transitive. Specifically, we will prove 
that the orbit of the point $(1/2, 0)$ is dense.

\begin{defi} Given a point $(\rho, \theta) \in M$ we define 
$$O_{f,m}
(\rho, \theta) = \lbrace (\rho,\theta), ..., f^{m-1} (\rho, 
\theta) \rbrace.$$
We recall that   a set $O \subset M$ is said to be $\epsilon$-dense in
$M$ if given any $y \in M$ there 
exists $x \in O$ such that $d(x,y) \leq \epsilon$.
\end{defi}

\subsubsection{The orbit of $(1/2, 0)$ under $f_n$.} Due to (\ref{yy1}) density of the orbit of $(1/2,0)$ under $f$ will follow if we prove that the set 
$O_{f_n, q_{n+1}} (1/2,0)$ is $1/n$-dense.

The following is evident if in our choice of $\a_{n+1}$ we 
require that $q_{n+1} \geq q_n^4$: For any $0 \leq p \leq 
q_n -1 $ and any $0 \leq k \leq q_n - 1$ and for $j=l$ and $j=r$ 
there exists
an integer $m(p,k,j) \leq q_{n+1}$ such that
$$\lbrace m \a_{n+1} \rbrace \in S_{p \over q_n} ( I_{n,k}^j).$$

Since $H_n (1/2,0) = (1/2,0)$ we have that
\begin{eqnarray*} f_n^{m(p,k,j) } (1/2,0) &\in&  H_n S_{p \over q_n} 
( I_{n,k}^j), \\
&=& H_{n-1} h_n S_{p \over q_n} (I_{n,k}^j) \\
&=& H_{n-1} S_{p \over q_n} h_n (I_{n,k}^j) \ \ \ \ {\rm from} \ 
(\ref{yyh2}) {\rm \ of \ Lemma \ } \ref{yyhn}, \\
&\subset& H_{n-1} S_{p \over q_n} (B_{n,k}^j)  \ \ \ \ \ \ \ {\rm from} 
\ (\ref{yyh4}) {\rm \ of \ Lemma \ } \ref{yyhn}.
\end{eqnarray*}

In conclusion, we proved that for any $0 \leq p \leq q_n -1 $ and any 
$0 \leq k \leq q_n - 1$, for $j=l$ and $j=r$  there exists
an integer $m(p,k,j) \leq q_{n+1}$ such that
\begin{eqnarray*} f_n^{m(p,k,j) } (1/2, 0) \in  H_{n-1} S_{p  \over 
q_n} (B_{n,k}^j). \end{eqnarray*}

Assuming that  ${||H_{n-1}||}_{C^1} \leq {q_n \over n}$ we have that 
diam$\left( H_{n-1} S_{p / q_n} (B_{n,k}^j)  \right) \leq 1 / n$ 
hence the partition ${\lbrace    H_{n-1} S_{p / q_n} (B_{n,k}^j) 
\rbrace}$, $0 \leq p \leq q_n -1 $, $0 \leq k \leq q_n - 1$, $j=l, 
r$, is a $1 / n$ grid of $M$.  Therefore the set $O_{f_n, q_{n+1}} 
(1/2,0)$ is indeed $(1/n)$-dense in $M$ and the proof of transitivity 
is complete. \carre

\subsection{Non ergodicity: $\mu ( {\phi}_f) = 0$.}

\subsubsection{Preliminary.} As a consequence of Property 
(\ref{yyh5}) of Lemma \ref{yyhn}  we get the following properties on 
the diffeomorphism $H_n = h_1 \circ ... \circ h_n$

\begin{coro} \label{yycor} For any $\rho_0 > 0$ there exists 
$n(\rho_0)$ such that for any $n \geq n(\rho_0)$ we have
\begin{eqnarray*}
H_n (D^r_{{1 \over 2} + \rho_0}) &\subset& H_{n(\rho_0)} (D^r_{{1 
\over 2} + {\rho_0 \over 2}}), \\
H_n (D^l_{{1 \over 2} - \rho_0}) &\subset& H_{n(\rho_0)} (D^l_{{1 
\over 2} - {\rho_0 \over 2}}).
\end{eqnarray*}
\end{coro}

\begin{proof}
Fix $\rho_0 >0$ and let $n(\rho_0) \geq 4$ be such that $\rho_0 \geq 
2 / n(\rho_0)$.


For $n \geq n(\rho_0)$ we have that $1 / n < \rho_0$ hence
by  Property (\ref{yyh5}) of Lemma \ref{yyhn}
\begin{eqnarray*}
h_n \left(D^r_{ {1 \over 2} + \rho_0} \right) \subset
  D^r_{ {1 \over 2} + \rho_0 - {1 \over 2^n}}, \end{eqnarray*}
and if $n-1 \geq n(\rho_0)$ we have $1 / (n-1) \leq 1 / n(\rho_0) < 
\rho_0 - 1 / 2^n$ hence

\begin{eqnarray*}
h_{n-1} h_n \left(D^r_{ {1 \over 2} + \rho_0} \right) \subset
  D^r_{ {1 \over 2} + \rho_0 - {1 \over 2^n} - {1 \over 2^{n-1}}}. 
\end{eqnarray*}
Since for any $l \geq n(\rho_0)$ we have $\rho_0 - {1 / 2^n} - ... - 
{1 / 2^l} \geq \rho_0 - 1/ 2^{n(\rho_0)-1} \geq \rho_0 / 2 \geq 1 / 
n(\rho_0)$ we can continue inductively and obtain that
\begin{eqnarray*}
h_{n(\rho_0)+1} \circ ... \circ h_n \left(D^r_{ {1 \over 2} + \rho_0} 
\right) \subset
  D^r_{{1 \over 2} + {\rho_0 \over 2} } \end{eqnarray*}
which implies the conclusion of the corollary, the case of $D^l_{ {1 
\over 2} - \rho_0}$ being similar.

\end{proof}

\subsubsection{Proof of $\mu (\phi_f) = 0$.}

We introduce the following set

$$E_l = \lbrace (\rho, \theta) \in M / \rho \neq {1\over 2} {\rm \ 
and \ for \ any \ } n \geq l+1, h_n(x)=x  \rbrace.$$

 From Property (\ref{yyh3}) in Lemma \ref{yyhn} we have that
\begin{eqnarray*} \mu(E_l) &\geq& 1 - {1 \over \sqrt{q_{l+1}}} - {1 \over 
\sqrt{q_{l+2}}} - ... \\
&\geq& 1 - {2 \over \sqrt{q_{l+1}}}.
\end{eqnarray*}
Since $H_l$ is measure preserving we also have $\mu(H_l(E_l)) \geq 1 
- 2 / \sqrt{q_{l+1}}.$  We will hence finish if we show that for any point 
$z \in H_l(E_l)$ the orbit of $z$ under $f$  is not dense in 
$M$. From the convergence condition (\ref{yy1})  in \S 
\ref{yyconvergence} it suffices to prove that there exists a fixed 
closed set $D(z)$ strictly included in $M$ such that for any $n$, 
$O_{f_n, q_{n+1}}(z) \in D(z).$

Fix $z \in H_l(E_l)$ and write  $z = H_l(z')$.  We can assume that 
for some $\rho_0>0$, $z' \in E_l \bigcap D_{1 / 2 + \rho_0}^r$, the 
other case $z' \in E_l \bigcap D_{1 / 2 - \rho_0}^l$ being similar.

For any $n \geq l$ and any $m \leq q_{n+1}$ we have that $f_n^m (z) = 
H_n S_{\a_{n+1}}^m (z').$ Since $S_{\a_{n+1}}^m (z') \in D^r_{1 / 2 + 
\rho_0 }$, Corollary \ref{yycor} implies that for any $n \geq 
n(\rho_0)$ we have $f_n^m(z) \in H_{n(\rho_0)} 
D^r_{1 / 2 + \rho_0 / 2} $. The latter being a fixed closed set 
strictly included in $M$ the proof is complete.   \carre

\section{Mixing flows on manifolds admitting a nontrivial ${\TT}^3$ action.}\label{sectionmixing}

All the constructions we have presented up to now share the property of being $C^{\infty}$ rigid, 
namely that the limit diffeomorphism $f$ satisfies $f^{q_n} \to\Id$  as $n\to\infty$ in  the $C^{\infty}$ topology for some sequence of integers $q_n$.
This is also the case for all the available constructions following  the approximation by 
conjugation method  introduced in \cite{AK}. It is an open problem whether there exists a circle action $\mS$  on a  compact manifold   and $\a$  such that the corresponding space $\A_{\a}$ contains a nonrigid  transformation. This would be a first step towards proving the existence of mixing maps in $\A_\a$ (Cf. Problem \ref{mixinginA}). 

\subsection{Time changes as limits of conjugacies}\label{timechangelimit}
The situation is different if we consider higher dimensional actions. Indeed, it is proved in  \cite{F1} that there are  linear flows on $\TT^3$  that can be made mixing via smooth (real analytic) reparametrization. Notice that reparametrizations of linear flows lie in the closure of conjugacies of periodic linear flows (and even of the standard flow $\mS ;\,\,S_t(x_1,x_2,x_3)=(x_1,x_2,x_3+t)\,\,\,(\mod 1)$
on ${\TT}^3$ by the following reasoning. 

Diophantine frequencies  are dense and  any smooth reparametrization of a Diophantine linear flow is  smoothly conjugated to a linear flow. Hence taking a very close rational approximation of the density  vector and applying the same conjugacy produces a periodic flow approximations the given one.
Moreover, every periodic linear flow on ${\TT}^n$ is conjugate to  the standard ``vertical'' flow 
$\mS$
along the last coordinate  on the torus albeit  with a constant  reparametrization to normalize the period to one.
The conjugacy is an automorphism of the torus and thus is not homotopic to identity.
Thus, every diffeomorphism from the mixing flow  obtained by a reparametrization  belongs to in the space $\A$   of limits 
of conjugates of elements of the action $\mS$.  

It is not clear however whether these on any other mixing maps in $\A$  belong to any of the spaces $\A_{\a}$. (Cf. discussion in Section \ref{closurespace}
and in particular Problem \ref{mixinginA}). Another open question more specific to the torus 
is the following.

\begin{prob}
Consider the  space $\tilde A$, the closure of conjugacies of the elements of the  vertical action $\mS$
on ${\TT}^n$  by volume preserving diffeomorphisms homotopic to identity.
Does the space $\tilde A$ contain a mixing diffeomorphism?
\end{prob}

\subsection{Formulation of results}
 In this section, given a manifold $M$ admitting a nontrivial smooth $\TT^3$ action preserving a smooth volume $\mu$ we construct $\mu$-mixing diffeomorphisms in the closure of conjugacies of elements of the action. To this end we combine the techniques of reparametrization of linear flows on $\TT^3$  with the explicit approximation by conjugation methods.  Besides working with higher dimensional actions and requiring a control for {\sl all iterates} of the limit diffeomorphism, the main difference with the former sections is that we need to work in restricted spaces $\A_{\a}$ with special Liouvillean vectors $\a \in \RR^3$ chosen {\sl a priori} such that the linear flow $R_{t \a}$ admits smooth mixing reparametrizations. Conditions of the type \ref{convergence} become more restrictive since we do not have complete freedom in the choice of the periodic approximations of $\a$ and explicit bounds on the conjugations used are required. This corresponds to quantitative  versions of  the approximation by  conjugation method similar to the approach of \cite{FKS}. 

The comments about mixing flows on the torus and the difference between the spaces  $\A$ and $\A_{\a}$ of course apply to this more general situation.

\begin{theo} \label{4 theorem} Assume $M$ is a differentiable manifold admitting a nontrivial smooth ${\TT}^3$ action ${\mathcal S} ={\lbrace S_v \rbrace }_{v \in \RR^3}$, $S_{v + k} = S_v$ if $k \in \ZZ^3$ and let $\mu$ be a smooth  volume preserved by ${\mathcal S}$. Then there exists a sequence $\gamma_n \in \QQ^3$ and a sequence $H_n$ of diffeomorphisms preserving $\mu$ such that the sequence $H_n S_{t\gamma_n} H_n^{-1}$ converges in the $C^\infty$ topology to a flow preserving $\mu$ and mixing for this measure. 
\end{theo}

\begin{rema}\label{remarkt2mixing} This is the optimal result for flows,  in the sense that a ${\TT^2}$ action is not sufficient:  for example, there are no mixing fixed point free flows on 
${\TT}^2$. One may expect that on a  manifold with a  nontrivial ${\TT}^2$ action  a discrete time version of  our construction should work. Namely, considering successive conjugations of some sequence of rational elements of the action it should be possible to obtain  mixing diffeomorphisms in the limit.  However,  at the moment we are not able to provide such a construction even on ${\TT}^2$. 
\end{rema}

\subsection{Reduction to $\TT^3 \times \DD^m$.}\label{4 reduction}$ $

Theorem \ref{4 theorem}
 will follow if we can prove it in
the particular case of a direct product  $ \TT^3 \times \DD^{m}$  with successive conjugacies $H_n$ having compact support in the interior of $\TT^3 \times  \DD^m$. 

This follows from the modification of Proposition \ref{reductiondirect} for the case of the torus action.

\begin{prop}Let $M$ be an $m+k$-dimensional 
differentiable manifold with an effective  action
  ${\mathcal T}$  of $\TT^k$ preserving a smooth volume.
$\mu$. Let $B$  be the union of the boundary of $M$  and the set of points with  a nontrivial stationary  stabilizer.  There exists a continuous surjective map $F: {\DD}^{m} \times \TT^k \to M$  with the following properties
\begin{enumerate}
\item The restriction of $F$ to the interior $ {\BB}^{m} \times \TT^k$ is a $C^{\infty}$ diffeomorphic embedding;
\item $\mu(F(\partial( {\DD}^{m} \times \TT^k))=0$;
\item $F(\partial( {\DD}^{m-1} \times S^1)\supset B$;
\item$F_*(\lambda)=\mu$;
\item $\mathcal T\circ F=F\circ\mS$, where $\mS$ is the standard vertical action on the direct product.
\end{enumerate}
\end{prop}

\subsection{Plan of the construction.} We will use the natural $\TT^3$ action on $\TT^3 \times \DD^m$, $S_v (\theta, r) = (\theta + v, r)$. We will not directly construct the mixing flow as a limit of conjugacies of periodic flows $H_n S_{t \gamma_n} H_n^{-1}$. Instead we will use successive conjugacies of a flow $T^t \times {\rm Id}_{\DD^m}$ where $T^t$ is the mixing flow on $\TT^3$ obtained in \cite{F1} via reparametrization of a Liouvillean linear flow on $\TT^3$. Theorem \ref{4 theorem} then follows since reparametrizations of linear flows are in the closure of conjugacies of periodic linear flows on ${\TT}^3$ as has been explained in Section \ref{timechangelimit}.

\subsection{Mixing reparametrizations of linear flows.} \label{4 notations} Since some specific  features of the mixing reparametrized flow constructed in \cite{F1} will be needed in the proof of Theorem \ref{4 theorem} we will summarize in this section the construction and mention the properties
that will be useful for the sequel.

\subsubsection{Notations and definitions.} 
For $\a \in \RR^3$, we denote by $R^t_\alpha$ the linear flow on the torus $\TT^{3}$  given by $ {d \theta} / {d t} = \alpha, $
where $\theta \in \TT^{3}$. For a continuous function $\phi : \TT^{3} \rightarrow \RR^{*}_{+}$ we
denote by  $R^t_{\alpha, \phi}$ the reparametrized flow  given by 
$  {d \theta} / {d t} = {\phi (\theta) \a}. $
 If the coordinates of $\alpha$ are rationally independent then the linear flow $R^t_{\alpha}$ is uniquely ergodic and so is the reparametrized flow that preserves the measure with density $1 / \phi$. In contrast, other properties of the linear flow such as mixing may be very sensitive to reparametrizations. The reparametrized flow $R^t_{\a, \phi}$ can be viewed as a special flow $T^t_{\a, \varphi}$ over $R_\a$ and under $\varphi$ where $\a$ this time is an $\RR^2$ vector and $\varphi$ is a strictly positive function over $\TT^2$. We denote by $M_{\a, \varphi}$ the space where the special flow acts which  is obtained as a quotient of $\TT^2 \times \RR$ by the equivalence relation $(x,s + \varphi(x)) \sim (x+\a,s)$ and that can be identified with $\TT^3$. Conversely any smooth special flow $T_{\a,\phi}$ can be viewed as a reparametrization of some linear flow on $\TT^3$. Henceforth, we will only work with special flows that are easier to manipulate than reparametrizations. The unique probability measure preserved by the special flow is given by the product of the Haar measure on the base $\TT^2$ and the Lebesgue measure on the fibers normalized by the mean value of $\varphi$. The $k$-dimensional Lebesgue measures will be denoted by $\lambda^{(k)}$. 

\subsubsection{Mixing special flows over a translation of $\TT^2$. }

Our choice of the translation $R_{\a,\a'}$ and the ceiling function $\varphi$ follows \cite{yoccoz1}
.  We take $(\a,\a') \in (\RR \setminus \QQ)^2$, such that
their sequences of best approximations $q_n$ and $q'_n$ satisfy for any $n \geq n_0(\a,\a')$
\begin{eqnarray}
q'_n \ \ &\geq&  e^{4q_n}, \label{5 sd 1} \\
q_{n+1} &\geq& e^{4q'_n}. \label{5 sd 2} \end{eqnarray}

\noindent The ceiling function we will consider is the following strictly positive real analytic function on ${\TT}^2$:             
\begin{eqnarray} \label{4 varphi}  \varphi(x,y)= 1+ {\rm Re} \left( \sum_{k=2}^{\infty} { e^{i2 \pi kx}
\over e^{k}} +
\sum_{k=2}^{\infty} {e^{i2 \pi ky} \over e^{k}} \right). \end{eqnarray}

In \cite{F1} it was proved that

\begin{theo} Under conditions (\ref{5 sd 1})--(\ref{4 varphi}), the special flow $T^t_{\a,\a',\varphi}$ is mixing for its unique invariant  probability measure.
\end{theo}

Since $\a, \a'$ and $\varphi$ are fixed, to alleviate the notations we will denote by $T^t$ the special flow $T^t_{\a,\a',\varphi}$ and by $M$ the space $M_{\a,\a', \varphi}$. We will need the following properties of  $T^t$, the first one of which is an obvious approximation of the flow by special flows over periodic translations and the other being a refinement of the mixing property that can be derived from the proof of mixing given in \cite{F1}.

  We will denote by $\a_n$ and $\a'_n$ the $n-$th convergents of $\a$ and $\a'$, i.e. $\a_n = p_{n} / q_{n}, \a'_n = p'_{n} / q'_{n}$. By $T^t_n$ we denote the special flow over the periodic translation  $R_{\a_n, \a'_n}$ and under the function $\varphi$.
Since $|\a - \a_{n}| \leq 1 / q_n q_{n+1}$ and  $|\a' - \a'_{n}|\leq 1 / q'_n q'_{n+1}$ we deduce

\begin{eqnarray} \label{5 S1} 
 ||T^t - T^t_n|| = O({t \over q_n q_{n+1}}). \end{eqnarray}

\subsubsection{The specialized mixing property of $T^t$.} \label{specspec} Let $\delta_n =
1 / q_n$. Consider a collection ${\mathcal C}_n = {\lbrace C_{n,j}
\rbrace}_{j \leq M_n}$ of disjoint sets in $ M$ of the form 
$$C_{n,j} = \left[ {i \over q_{n}} , {i+ 1 \over q_{n}} \right] \times
\left[ {i' \over q'_{n}} , {i'+ 1 \over q'_{n}} \right] \times  \left[ s,
s+ \delta_n \right]$$  
for some $i=i_j \leq q_{n} -1$ and $i'=i'_j \leq q'_{n} -1$ and $s=s_j$
satisfying
$$ s+ \delta_n \leq \min_{(x,x') \in [i_j / q_{n} , (i_j+ 1) / q_{n}]
\times
[i'_j / q'_{n} , (i'_j + 1) / q'_{n}]} \varphi(x,x').$$
We  take the collection ${\mathcal C}_n$   such  that $\mu( \bigcup
C_{n,j}) \geq 1 - \epsilon_n$ with $\epsilon_n = O (1/q_n)$.

\begin{defi} An {\sl $\epsilon$-grid of a  set $(A, \mu)$}  is any
partition  of $A$  into disjoint measurable sets of diameter less than
$\epsilon$. An {\sl $\epsilon$-grid of a collection of disjoint similar
sets ${\lbrace C_{j} \rbrace}_{j \leq M}$} is given by some
$\epsilon$-grid of $C_1$ identically repeated inside each $C_j$, $j \leq
M$. 
\end{defi}

\begin{defi} \label{5 even} For a pair of strictly positive numbers
$(\epsilon, \eta)$, a set $B \subset A$ is said to be {\sl $(\epsilon,
\eta)$-uniformly-distributed inside $A$} if there exists an
$\epsilon$-grid of $A$ such that for each atom $A_k$ of the grid we have 
$$\left| {\mu(B \bigcap A_k)  \over \mu(B) } - {\mu(A_k) \over \mu(A)}
\right| \leq \eta  {\mu(A_k) \over \mu(A)}.$$

 
\end{defi}

\begin{defi} \label{5 uniform} For a pair of strictly positive numbers
$(u,v)$, a set $B$ is said to be {\sl $(u,v)$-identically-distributed in a
collection  of disjoint similar sets ${\lbrace C_j \rbrace}_{j \leq M}$}
if there exists a $u$-grid of the collection   ${\lbrace C_{j} \rbrace}_{j
\leq M}$ such that for every pair $j,j' \leq M$ and for any $k \leq K$
($K$ referring to the number of atoms inside each $C_j$), we have  
$$ \left| \mu(B \bigcap C_{j,k}) -  \mu(B \bigcap C_{j',k}) \right| \leq v
\mu(B \bigcap C_{j,k}).$$
\end{defi}

\begin{rema} In the above definition, the measurable sets $B$,  $A$,
$A_k$, $C_{j,k}$ might be of different dimensions and $\mu$ should be
replaced when necessary by its appropriate conditionals.
\end{rema}

\begin{rema} \label{rr33} A set $B$ that is
$(u,v)$-identically-distributed inside the collection  ${\mathcal C}_n=
{\lbrace C_{n,j} \rbrace}_{j \leq M_n}$ is $(2/q_n,
v)$-uniformly-distributed in $M$ (this is due to the fact that the
diameter of the sets $C_{n,j}$ is less than $2 / q_n$). Moreover, because
each set $C_{n-1,j}$ is almost a union of sets of ${\mathcal C}_n$ then $B
\bigcap C_{n-1,j}$ is $(2/q_n, 2v)$-uniformly distributed inside
$C_{n-1,j}$.
\end{rema}

The specific mixing property of $T^t$ that will be crucial in the sequel
is the following

\begin{prop} \label{5 S2} Let $\xi_n$ be a set in $M_{\a, \varphi}$ of the
form $I \times J \times \lbrace s \rbrace $ where $I$ and $J$ are
intervals of size $1 / e^{q'_n}$. Then for any $t \geq e^{3q'_n}$ we have
that $T^t \xi_n$ is $(1/e^{q'_n}, 1/ n^2)$-identically-distributed in the
collection ${\mathcal C}_n= {\lbrace C_{n,j} \rbrace}_{j \leq M_n}$.
\end{prop}

\subsection{The successive conjugations scheme.}

We will use the same notation $T^t$ and $T^t_n$  to indicate the special
flows over $T^t_{\a,\a', \varphi}$ and $T^t_{\a_n, \a_n', \varphi}$ as
well as the flows $T^t \times {\rm Id}_{\DD^m}$ and  $T^t_n \times {\rm
Id}_{\DD^m}$. 
 For $r \in (0,1]$, we denote by $D(0,r)$ the closed disc of radius $r$
inside $\DD^m$. The general scheme of the construction is the
following:  We construct volume preserving diffeomorphisms $h_n$ compactly
supported on  $M_{\a_n, \a_n',\varphi} \times \DD^m$ of the form 
\begin{eqnarray} \label{4gn} h_n(x,x',s,r) = (x,x',s, g_{n,x,x'}
(r)) \end{eqnarray} 
where for each $(x,x') \in \TT^2$ $g_{n,x,x'}$ is a compactly supported
volume preserving diffeomorphism of $\DD^m$ satisfying 
\begin{eqnarray} \label{44gn} g_{n,x,x'}(r) = g_{n, x + \a_{n},
x'+\a'_{n}} (r) \end{eqnarray}
for all $(x,x',r)$. This  implies that $h_n$ commutes with $T^t_{n}$. We
can choose $g_{n,x,x'}$ increasingly sensitive with respect to $x$ and
$x'$ but in a controlled as follows 
\begin{prop} \label{5 hn} It is possible to construct the  diffeomorphism
$h_n$ as in (\ref{4gn})--(\ref{44gn}) such that the following holds

\begin{enumerate}

\item \label{5 h0} $h_n$ is volume preserving and compactly supported,

\item \label{5 h1} ${h}_n T^t_{n} = T^t_{n} h_n$,


\item  \label{5 h2} For any $l \in \NN$, for $n$ large enough
${||{h}_n||}_{C^l}  \leq e^{{1 \over 2} \sqrt{q'_{n}}}$.

\item  \label{5 h3} For any $r \in D(0, 1 - {e^{-n}}) \subset \DD^m$, for
any $j \leq M_{n}$ we have  that \\
${h}_n \left(C_{n,j} \times \lbrace r \rbrace  \right) $ is
$({e^{-q'_{n-1}}}, {1 \over n})$-uniformly-distributed inside $
C_{n,j} \times \DD^m$.
  in the sense that there exists a partition of $\DD^m$ in sets having
diameter less than ${e^{-q'_{n-1}}}$ such that for any atom $D$ of the
partition we have  
$$\left| \lambda^{(3)} \left({h}_n \left(C_{n,j} \times \lbrace z \rbrace
\right) \bigcap  (C_{n,j} \times D) \right) - \lambda{(3)}
(C_{n,j}) \lambda^{(m)}(D) \right| \leq {1 \over n} \lambda^{(3)}
(C_{n,j}) \lambda^{(m)} (D).$$

\end{enumerate}

\end{prop}

We now consider the flows $T^t$, $T^t_n$ and the diffeomorphisms $h_n$ as
flows and diffeomorphisms  on $M \times \DD^m$ and define $H_n = h_1 \circ
... \circ h_n$ and consider the sequence $H_n T^t H_n^{-1}$.

 From the above proposition  we can easily deduce the following properties
of $H_n$

\begin{prop} \label{5 Hn} The diffeomorphisms $H_n$ are compactly
supported in the interior of $M \times \DD^m$ and satisfy

\begin{enumerate}

\item \label{5 H0} $H_n$ is volume preserving,

\item \label{5 H1} For any $l \in \NN$, for $n$ large enough
${||{H}_n||}_{C^l} \leq e^{\sqrt{q'_{n}}}$,


\item  \label{5 H2} For any $z \in D(0, 1 - e^{-n}) \subset \DD^m$, for
any $j \leq M_{n}$ we have  that \\
${H}_n \left(C_{n,j} \times \lbrace z \rbrace  \right) $ is  $({e^{-{1
\over 2} q'_{n-1}}}, {1 / n})$-uniformly-distributed inside $
C_{n,j} \times \DD^m$.
 
\end{enumerate}

\end{prop}

\noindent {\sc Proof.} Property (\ref{5 H1}) follows from (\ref{5 h2}) of
Proposition  \ref{5 hn}. Property (\ref{5 H2}) then follows from (\ref{5
h3}) of Proposition  \ref{5 hn} applied to $h_n$ and from (\ref{5 H0}) and
(\ref{5 H1}) of Proposition \ref{5 Hn} applied to $H_{n-1}$. \carre

\subsection{Convergence and approximation.} 

\begin{prop} \label{5 approximation} The sequence of flows  $H_n T^t
H_n^{-1}$ converges in the $C^\infty$ topology to a volume preserving flow
$\overline{T}^t$ of $M \times {\DD}^m$. Moreover, for $t \in [e^{3q'_n},
e^{3q'_{n+1}}]$ and $l \in \NN$ we have for $n$ large enough
$$ {||H_n T^t H_n^{-1} - \overline{T}^t ||}_{C^l} \leq {1 \over 2^n}.$$
\end{prop}

\noindent {\sc Proof.} From (\ref{5 h1}) of Proposition \ref{5 hn} we have
that
$$H_{n+1} T^t_{ {n+1}} H_{n+1}^{-1} = H_n T^t_{ {n+1}} H_n^{-1},$$
hence 

$$ {||H_n T^t H_n^{-1} - H_{n+1} T^t H_{n+1}^{-1} ||}_{C^l} \leq \psi_n +
\psi'_n$$
where 
\begin{eqnarray*} \psi_n &=& {||H_n T^t H_n^{-1} - H_{n} T_{n+1}^t
H_{n}^{-1} ||}_{C^l}, \\ \psi_n' &=& {||H_{n+1} T_{n+1}^t H_{n+1}^{-1} -
H_{n+1} T^t H_{n+1}^{-1} ||}_{C^l}, \end{eqnarray*}
hence from (\ref{5 S1}) and (\ref{5 H1}) of Proposition \ref{5 Hn} we
conclude that 
$$ {||H_n T^t H_n^{-1} - H_{n+1} T^t H_{n+1}^{-1} ||}_{C^l} = o \left( {
e^{{q'_{n+1}}} t \over q_{n+1} q_{n+2}} \right).$$
Since $q_{n+2} \geq e^{4q'_{n+1}}$ the proposition follows. \carre

\subsection{Proof of mixing.} It remains to prove the following 

\begin{prop} The volume preserving flow of $M \times {\DD}^m$
$$\overline{T}^t = \lim_{n \rightarrow \infty} H_n T^t  H_n^{-1}$$
 is mixing.
\end{prop}

\subsubsection{} \label{5 mix} The following is an almost straightforward Fubini argument  
\begin{lemm}   \label{5 fub}

and $R$ are open balls in $\DD^m$ and ${\TT}^2$  and $\delta \leq c$,

If there is a sequence of measurable collections of subsets of  $M \times
\DD^m$, $\lbrace {\mathcal P}_t \rbrace$
  satisfying the following 
$$\lbrace {\mathcal P}_t \rbrace   \mathop{\longrightarrow}\limits_{ t
\rightarrow \infty }
  \varepsilon \ \ ({\rm \ partition \ into \
points  \ }), $$
and such that for any open set $A \subset  M \times \DD^m$ and any
$\varepsilon > 0$, there exists $t_0$ such that for $t > t_0$ we have for
every atom $P_{t} \in \lbrace {\mathcal P}_t \rbrace$ 

\begin{eqnarray*}  \left| \mu \left( P_t \bigcap
\overline{T}^{-t}(A)  \right)  -  \mu(P_t) \mu(A) \right| \leq \epsilon
\mu(P_t), \end{eqnarray*}
then the flow $\overline{T}^t $ is mixing.

\end{lemm}

\begin{rema} The atoms of the collection $\lbrace {\mathcal P}_t \rbrace $
might have positive codimension and the Lemma remains valid if appropriate
conditional measures are used in (\ref{5 mix}).
\end{rema}

\subsubsection{\sl Choice of the collection  ${\mathcal P}_t$.} \label{5 boundary}  Let $t \in
\RR_+$ and take $n$ such that  $t \in [ e^{3q'_{n}}, e^{3q'_{
{n+1}}}]$. We first consider a collection ${\mathcal U}_n$  converging to
the point partition in $M_{\a,\varphi} \times \DD^3$  consisting of sets
of the form $\xi \times  \lbrace s, r \rbrace $ where $\xi$ is a square on
${\TT}^2$ with side $1 / e^{q'_n}$. Since the diffeomorphism $h_n T^t$ is
volume preserving we remain with a collection converging to the point
partition if we discard the elements  
$\xi \times \lbrace s,r \rbrace $
that do not satisfy 
\begin{eqnarray*}  \lambda^{(2)} \left[ (\xi \times
\lbrace s,r \rbrace \bigcap T^{-t} h_n^{-1} ( M_{\a,\varphi} \times D(0,
1- e^{-(n-1)}))  \right] \geq (1 - {1 \over n})  \lambda^{(2)}
(\xi). \end{eqnarray*}

Finally, we consider the collection ${\mathcal P}_t$ whose atoms are the
sets $H_n(\xi \times \lbrace s, r \rbrace )$ where $ \xi \times \lbrace
s,z \rbrace$ is as above. Due to the fact that ${||H_n||}_{C^1} \leq
e^{\sqrt{q'_n}}$ we have that $ {\mathcal P}_t$  converges to the
partition into points of $M_{\a,\varphi} \times \DD^m$ as $t$ goes to
infinity.

\subsubsection{} \label{5 mixn} Fix an open set  $A \subset M \times \DD^m$ and $\varepsilon > 0$. We
finish if we check (\ref{5 mix}) of Lemma \ref{5 fub} for the elements of
${\mathcal P}_t$ as $t \in [e^{3q'_n}, e^{3q'_{n+1}}]$ goes to infinity or
equivalently as $n \rightarrow \infty$. To have (\ref{5 mix}) for this
range of $t$, we are reduced in light of Proposition \ref{5 approximation}
to proving that for $n$ large enough 
\begin{eqnarray*}  \left| \lambda^{(2)} \left( \xi\times
\lbrace s,r \rbrace    \bigcap  T^{-t} H_n^{-1}  (A) \right)   -
\lambda^{(2)} (\xi) \mu(A) \right| \leq \varepsilon \lambda^{(2)}
(\xi) \end{eqnarray*} 
where $\xi$ is a square in $\TT^2$ of side $1/ e^{q'_n}$ and $\xi \times
\lbrace s,r \rbrace$ satisfies (\ref{5 boundary}). 

\subsubsection{} \label{55 mixn} Since the collection ${\mathcal C}_{n-1} = {\lbrace C_{n-1,j} \rbrace }_{j
\leq M_{n-1}}$ converges to the point partition, (\ref{5 mixn})  will
follow if for an arbitrary measurable set  $\Delta \subset \DD^m$ we prove that
for $n$ large enough we have for any $j \leq M_{n-1}$   
\begin{eqnarray*}  \left| \lambda^{(2)} \left(\xi\times \lbrace s,r \rbrace
\bigcap     T^{-t} H_n^{-1} (C_{n-1,j} \times \Delta)  \right)  -
\lambda^{(2)} (\xi) \lambda^{(3)}(C_{n-1,j}) \lambda^{(m)} (\Delta) \right| \\ 
 \leq \varepsilon \lambda^{(2)}
(\xi) \lambda^{(3)}(C_{n-1,j}).  \hspace{1.5cm} \end{eqnarray*}
 
\subsubsection{} \label{5reduction} Having fixed $t$ and the atom $\xi\times \lbrace s,r \rbrace \in {\mathcal
P}_t$ we set
$$U := h_n \left(T^t (\xi\times \lbrace s,r \rbrace) \right)$$
and denote by $\tilde{\lambda}^{(2)}$ the measure  ${(h_n T^t)}_{*} (
\lambda^{(2)}$.

The equation we have to establish (\ref{55 mixn}) becomes 
 \begin{eqnarray*}  \left| \tilde{\lambda}^{(2)} \left(  U \bigcap
H_{n-1}^{-1}  (C_{n-1,j} \times \Delta)  \right)  -  \tilde{\lambda}^{(2)}
(U) \lambda^{(3)}(C_{n-1,j}) \lambda^{(m)} (\Delta) \right| \\
 \leq \varepsilon \tilde{\lambda}^{(2)}
(U) \lambda^{(3)}(C_{n-1,j}).  \hspace{1.5cm} \end{eqnarray*}

\subsubsection{} \label{clam} From Proposition \ref{5 S2} we recall that since $t \geq e^{3q'_{n}}$, then   the set $T^t (\xi \times \lbrace s
\rbrace)$  is  $(1/e^{q'_n}, 1/n^2)$-identically-distributed
in ${\mathcal C}_n$.
This  means that  $T^t (\xi \times \lbrace s
\rbrace)$  has essentially the same trace inside each $C_{n,j}$ with precision (or resolution) $e^{-q'_n}$. The diffeomorphism $h_n$ introduced in Proposition \ref{5 hn} has an important deformation effect on the scale of the sets $C_{n,j}$ (Property \ref{5 h3}) but negligible on the scale $e^{-q'_n}$ (Property \ref{5 h2} asserting that ${\parallel h_n \parallel}_{C^1} \leq e^{\sqrt{q'_n}}$). On the other hand by   (\ref{4gn}) and (\ref{44gn}) we have that $h_n (C_{n,j} \times \DD^m) = C_{n,j} \times \DD^m$  for any $j,j' \leq M_n$ and that 
${h_n}_{|C_{n,j} \times \DD^m} \equiv {h_n}_{|C_{n,j'} \times \DD^m}$. We obtain therefore the following 

\noindent {\sl Claim.} There exists a partition of $\DD^m$ 
 in open sets of diameter less than $1 / e^{{1 \over 2}q'_n}$, $\DD^m =
\bigcup_{l=1}^{l_n} D_l$, such that for each $l \leq l_n$, the set $U_l :=
U \bigcap (M \times D_l)$  is $(1/e^{ {1\over2} q'_n},
1/n)$-identically-distributed   in   the  collection   ${\lbrace C_{n,j}
\times D_l \rbrace}_{j \leq M_n} \subset M \times \DD^m.$

\subsubsection{} \label{52222} Since the collection ${\mathcal C}_n$ is almost a partition of $M$ that is much finer than ${\mathcal C}_{n-1}$ (Cf. \S \ref{specspec}), it follows from the claim above that 
for any $j \leq M_{n-1}$ 
\begin{eqnarray*} \left| 
\tilde{\lambda}^{(2)} \left( U_l \bigcap (C_{n-1,j} \times D_l)  \right) -
\tilde{\lambda}^{(2)}(U_l) \lambda^{(3)}(C_{n-1,j}) \right| \\  \label{ul}
\leq  {3 \over n}
\tilde{\lambda}^{(2)}(U_l) \lambda^{(3)}(C_{n-1,j}) \hspace{1.5cm} 
\end{eqnarray*}

\subsubsection{} \label{layl}  From Remark \ref{rr33} the Claim \ref{clam}  also implies that the set $U_l$ is $(2 /
q_n, 1 / n)$-uniformly-distributed in $M \times D_l$
and more specifically that for any $j \leq M_{n-1}$  
the  set  
\begin{eqnarray*}  U_l \bigcap (C_{n-1,j} \times D_l)  {\rm \ is \ }  (2 / q_n, 2 /
n)- {\rm uniformly-distributed \ in \ }   C_{n-1,j} \times D_l.  \end{eqnarray*}

\subsubsection{} From our choice of the collection ${\mathcal P}_t$ in (\ref{5 boundary}) we have that 
\begin{eqnarray} U \sim  \bigcup_{D_l \subset D(0, 1 - e^{-(n-1)}) } U_l
\label{5number}
\end{eqnarray} 
in the sense that the ratio between the measures of the left hand side set
and the right hand side set converges to 1 as $n$ goes to $\infty$. So we
assume from now on that $l$ is such that $ D_l \subset D(0, 1 - e^{-(n-1)})$.
We fix a measurable set $\Delta \in \DD^m$ and $j \leq M_{n-1}$ and proceed to estimate $\tilde{\lambda}^{(2)} (U_l \cap H_{n-1}^{-1} (C_{n-1,j} \times \Delta))$ in the scope of proving (\ref{5reduction}).

\subsubsection{} \label{678} Since $h_n (C_{n,j} \times \DD^m) = C_{n,j} \times \DD^m$  (from (\ref{4gn})) and  $U_l = U \cap (M \times D_l)$ we have that 
$$U_l \cap H_{n-1}^{-1} (C_{n-1,j} \times \Delta)   =   \left( U_l \cap 
(C_{n-1,j} \times D_l) \right)  \bigcap H_{n-1}^{-1} (C_{n-1,j} \times \Delta).$$

\subsubsection{} \label{f3} Because $D_l  \subset D(0, 1 - e^{-(n-1)})$, Property
(\ref{5 H2}) in Proposition \ref{5 Hn} asserts that $H_{n-1} (C_{n-1,j} \times D_l)$ is $(e^{-{1 \over 2} q'_{n-1}}, 1/(n-1))$-uniformly-distributed inside $C_{n-1,j} \times \DD^m$. From (\ref{layl}) and the fact that ${||H_{n-1}||}_{C^1} \leq e^{\sqrt{q'_{n-1}}}$ we deduce that $H_{n-1}  ( U_l \cap 
(C_{n-1,j} \times D_l) )$ is   $(e^{-{1 \over 4} q'_{n-1}}, 4/n)$-uniformly-distributed inside $C_{n-1,j} \times \DD^m$ and get with (\ref{678})

\begin{eqnarray*} 
\tilde{\lambda}^{(2)}  \left[ U_l \cap H_{n-1}^{-1} (C_{n-1,j} \times \Delta)  \right] \sim  \tilde{\lambda}^{(2)}  \left(U_l \bigcap (C_{n-1,j}
\times D_l) \right) \lambda^{(m)}(\Delta) \end{eqnarray*}
which with (\ref{52222}) becomes 
\begin{eqnarray*} 
\tilde{\lambda}^{(2)}  \left[ U_l \cap H_{n-1}^{-1} (C_{n-1,j} \times \Delta)  \right] \sim 
\tilde{\lambda}^{(2)}(U_l) \lambda^{(3)} (C_{n-1,j})  \lambda^{(m)}(\Delta). 
\end{eqnarray*}
 
\subsubsection{} Having 
 (\ref{5number}) in mind we can sum in (\ref{f3}) over $l$
such that $D_l \subset D(0, 1 - e^{-(n-1)}) $ and get 

\begin{eqnarray*} \tilde{\lambda}^{(2)} \left( H_{n-1} U \bigcap
(C_{n-1,j} \times \Delta) \right) 
&\sim&
\tilde{\lambda}^{(2)}(U) \lambda^{(3)} (C_{n-1,j}) \lambda^{(m)}(\Delta). \end{eqnarray*}
This gives equation (\ref{5reduction})  and the proof of mixing is hence
accomplished.   \carre

\section{Other directions and more open problems}\label{lastsection}
We discuss two areas where  some of the most interesting and challenging problems related to the approximation by conjugation  method appear: analytic constructions
and the problem of smooth realization. 

\subsection{Real--analytic constructions}The approximation by conjugation me-thod is a major  source of constructing  examples of smooth dynamical systems with prescribed ergodic or topological properties; other methods are  skew products and time--changes. In the latter two cases, and except for some special exotic constructions \cite{F5}, there is usually not much difference between the cases of
sufficiently high finite differentiability and  real analyticity. Constructions in \cite{furstenberg},  \cite[Sections 12.3, 12.4, 13.3, 13.4]{K1}, \cite{F1, F2, F3, F4}  provide typical examples. In a somewhat crude way  this can be explained  by essential linearity of the time change and skew product constructions.  The conjugation by approximation construction is essentially nonlinear and it is based  on convergence of  maps  obtained 
from  certain standard maps by  wildly diverging conjugacies. 

Here a great difference  between the differentiable and real analytic
maps becomes apparent.   In the former case as long as the  conjugating diffeomorphisms  are sufficiently smooth albeit with huge derivatives  fast  convergence  of the approximations is achieved  by conditions like  
\eqref{convergence}. Since  no control over $C^r$ norms is required it is usually easy 
to construct $C^{\infty}$ diffeomorphisms creating an approximate picture  required.
Of course usually smooth maps can be approximated by real analytic ones. Such maps if they 
have very large derivatives in the real domain 
will usually have singularities in a fairly small complex neighborhood of it.
But if either  a map $h$  or its inverse has complex singularities close to the real domain and  $S_t,\,\,0\le t\le t_0, \,\, S_0=\Id$ is 
an analytic family  then the  family $h\circ S_t\circ h^{-1}$ is expected to  have singularities
close to the real domain  for {\em any} $t > 0$.  Thus the domain of    analyticity  for maps $f_n$ in \eqref{conjugations} will shrink
considerably at any step of the constructions and the limit map $f$ will not be analytic.

There are  two  possible strategies for  overcoming this problem in order to make a version  of  the approximation by conjugation  method work in the real--analytic category.
\begin{enumerate}
\item Find conjugacies  $h_n$  which are analytic  with their  inverses in a  uniform complex neighborhood  of the real domain
\item Compensate  singularities of the conjugacies  $h_n$ and their inverses   in such a way 
that the  maps $f_n$ have a  much larger domain of analyticity than the  conjugacies.
\end{enumerate}

\subsubsection{Uniquely ergodic diffeomorphisms of spheres} The first method has a good chance of working when there is  sufficiently large supply
 of maps commuting with the action $\mS$ which can be modified in a controlled way
 maintaining analyticity in the large domain. This happens   for skew  products  where the $S^1$ acts in the fibers  and in the time-changes of linear flows on the torus. A less obvious case of this situation  was found by the second author \cite{K3}. We describe one of the results  of that work in a special case.
 
Consider the standard embedding of the sphere $S^{2n-1}$ into $\RR^{2n}$
and  the standard complexification  $\RR^{2n}\subset\CC^{2n}$.
The vector-field   defined in Euclidean coordinates as 
$v_0(x_1,\dots,x_{2n})=2\pi(x_2,-x_1,\dots,x_{2n},-x_{2n-1})$  defines a  
linear action $\Phi=\{\varphi ^t,\,\,t\in\RR,\,\, \varphi^1=\text{Id}\}$ of 
the circle
$S^1$ 
In Euclidean coordinates 
$$
\aligned
\varphi^t(x_1,\dots, x_{2n})= &(\cos2\pi tx_1+\sin2\pi tx_2, -\sin2\pi tx_1+\cos2\pi tx_2,\dots,
\\
&  \cos2\pi tx_{2n-1}+\sin2\pi tx_{2n}, -\sin2\pi tx_{2n-1}+\cos2\pi tx_{2n}).
\endaligned
$$ 
We will use the same notations $v_0$ and $\varphi^t$  for extensions to 
$\CC^{2n}$
or its subsets.

We will call a function on $S^{2n-1}$  {\em entire}  if it extends  to a holomorphic function
in $\CC^{2n}$.  A map  $f:S^{2n-1}\to S^{2n-1}$ is called  {\em entire} if 
its coordinate functions are entire.  A diffeomorphism  $f:S^{2n-1}\to S^{2n-1}$ is   entire if both $f$ and $f^{-1}$ are entire maps. Invertible  linear maps are obviously entire diffeomorphisms.
Notice that product of entire diffeomorphisms is an entire diffeomorphism.

\begin{theo}\label{theoremmainsphere} 
For any $t_0,\,\,0\le t_0\le1,\,\,\epsilon>0$ and any compact set 
$K\subset\CC^n$ there
exists an entire  diffeomorphism  $f:S^{2n-1}\to S^{2n-1}$, preserving 
Lebesgue measure and uniquely ergodic whose extension to $\CC^{2n}$
satisfies
\begin{equation}\label{eqclosenessCn}
\max \{ \max_{z\in K}|f(z)-\varphi^{t_0}(z)|, \max_{z\in K}|f^{-1}(z)-\varphi^{-t_0}(z)|\}\le\epsilon.
\end{equation}
\end{theo}

\begin{rema}Notice that while in the $C^{\infty}$ category  the  approximation by conjugation   method  as outlined in Section \ref{general scheme} allows to produce a great variety of perturbations of a rotation $\varphi_t$  with  interesting ergodic and topological properties, even minimal real analytic diffeomorphisms of $S^3$ were not known before. 
\end{rema}

This result extends to compact Lie groups and some of their homogeneous spaces with an 
action of $S^1$ by left translations. 

Let us sketch some of the ideas of the proof of Theorem \ref{theoremmainsphere} in the case of $S^3$.
There is a natural identification  of $S^3$ with the group $SU(2)$  which  identifies  $\Phi$ with an action by  left translations. The space of orbits is naturally identified with $S^2$ (the Hopf fibration).
Right translations   on $SU(2)$ act on  the factor space  as rotations with respect to
a  natural  metric which comes from a bi-invariant metric on $SU(2)$. 
Neither an individual translation  nor a one--parameter group is transitive  
but  of course  the whole right action is transitive.  Moreover,  every 
right translation also extends to a linear map of  $\CC^4$ and is hence 
entire. 

At $n$th iterative step of the construction one produces an entire   
diffeomorphism $h_n$ of $S^3=SU(2)$  which  commutes with $\varphi_{1/q_n}$. To imitate   the uniform distribution in the factor space one uses  twisted right translations. First, one  fixes a one--parameter subgroup $r_t$ of  right translations of $SU(2)$, or equivalently, an axis of rotation on $S^2$. On each orbit  $\O$ of this   one--parameter group  the map $h_n$ is  a  certain translation $r_{f(\O)}$
where the ``twisting  function '' $f$ is entire and is invariant under the left translation $\varphi_{1/q_n}$.
By making this function  sufficiently uniformly distributed  $\mod 1$ one achieves a distribution of  orbits of the  conjugate of $\Phi$ (and hence of  any  of its sufficiently  large finite subgroups)
which projects to $S^2$ into a distribution close to the  distribution  uniform on  orbits of 
$r_t$. This distribution however  is far from uniform. Further approximation is achieved
by taking another one--parameter subgroup and averaging the distribution  along its  orbits.
In order to implement  this one takes a number $q_n^{'}>>q_n$ depending on the geometry of the twisting function and constructs another twisting function along the orbits of the new subgroup
and  invariant with respect to $\varphi_{1/q_n^{'}}$, taking the conjugacy as the  product of two 
twisted right translations. This product is still entire and it takes  the  distribution of  orbits  of $\Phi$
and hence any of its sufficiently large finite subgroups closer to uniform. In fact, in order to achieve
almost uniform distribution of {\em every} orbit one need to either take more than two  steps within the  inductive step, or  choose  one--parameters subgroups of right translations  at different inductive steps carefully.

\begin{rema}This construction produces  conjugacies which are far from  preserving  natural fundamental domains of $\varphi_{1/q_n}$. In other words, the choice of $h_n$ is not special in the sense of Section \ref{sbsspecial} (See also Section \ref{almostallall}). This makes controlling subtle  ergodic properties  of the resulting uniquely ergodic diffeormorphism problematic.
\end{rema}
 
 Here is a typical open problem related with this difficulty.
 
 \begin{prob} Does there exist a real analytic volume preserving  diffeomorphism of $S^{2n-1},\,\,n\ge 2$ which is measurably conjugate to an irrational rotation of the circle?
 \end{prob}
 
\subsubsection{Exotic holomorphic maps in one complex variable}
The second method requires  good understanding of the  geometry of  special holomorphic maps
in the complex domain
since  just approximating  the conjugacies  by   some standard regular  objects  even if successful
usually produces   singularities close to the real domain  for inverse maps. In one complex variable there is a sufficient understanding  of such special geometry. This allowed 
R. Perez-Marco and   J.-C. Yoccoz \cite{P1, P2, P3, Y2},  developing ideas   which originally appeared in   Yoccoz'  thesis) to apply  approximation by conjugation  method
and construct analytic circle  diffeomorphisms with  Liouvillean  rotation numbers 
and exotic properties (e.g. $C^{\infty}$ but not analytic  conjugacy to a rotation).

\subsubsection{Other situations}  Beyond the two  cases  described above the 
problem of applicability of  the approximation by conjugation scheme of
Section \ref{general scheme} remains open.
 This 
is the case  
already in the setting of  volume preserving diffeomorphisms in dimension two
where the majority of examples described in this paper appear.
Certain characteristic open problems which appeared already at the time  
of writing \cite{AK} can be summarized as follows.

\begin{prob} Does there exist an area preserving topologically transitive  real analytic diffeomorphism  $f$ of the  disc $\DD^2$ with  either of the following properties:
\begin{enumerate}
\item the restriction of $f$ to the boundary  is an irrational rotation;
\item $f$ has zero topological entropy;
\item $f$  is $C^2$ close to the identity?
\end{enumerate}
\end{prob}

We finish our discussion of analytic constructions with   several comments concerning  the above problem.

In the case (1) the rotation number must be Liouvillean; otherwise there are invariant circles
as was explained in Section \ref{Discdiffeo}. 

Applicability  of the approximation by conjugation method in the analytic setting in the  case of manifolds with  boundary such as $\DD^2$  depends on our ability to  construct real analytic diffeomorphisms of such manifolds with very large derivatives  which can be extended 
together with their inverses  to  fixed complex neighborhoods. Here is a characteristic question of this kind. 

\begin{prob}
Does there exist an $r>0$ such that for any point $p=(x,y)\in\text{Int} \,\DD^2$ there exists a real analytic area preserving diffeomorphism $f_p:\DD^2\to\DD^2$  which extends together with its inverse to  $r$-neighborhood of $\DD^2$ in the complex domain and such that $f_p(0,0)=(x,y)$?
\end{prob}

Notice that for $p$ close to the boundary such a diffeomorphism $f_p$ must have  very large derivatives. 

Known examples of real analytic area preserving ergodic (or just topologically transitive)
diffeomorphisms of $\DD^2$ are based on modifications on the $C^{\infty}$ construction first introduced in \cite{K}; see \cite{G} and \cite{ LS}. Dynamically these examples are non-uniformly hyperbolic,  they have  positive topological entropy, and with respect to the smooth measure they are Bernoulli  (or direct  products of Bernoulli with finite permutations).  On the boundary  these diffeomorphisms
are either identity or  have  rational rotation number with finitely many periodic orbits.  Such examples can be made $C^1$ close to  identity (first  observed by Michel Herman; private communication)  but  constructions of this kind  are not compatible with $C^2$ closeness to identity.

\subsection{The problem of smooth realization}\label{sectionsmooth}
In its most basic and at the same time general form  the smooth realization problem may be formulated as follows;

\begin{prob}[Smooth realization] Given a measure preserving transformation $T$  of a Lebesgue space $(X,\mu)$ does there exist a ($C^{\infty}$) diffeomorphism  $f$ of a compact manifold $M$
preserving  a smooth volume $\nu$ such that $(f,\nu)$ is measurably isomorphic  to $(T, \mu)$?
\end{prob}

There are numerous variations and specializations of the  problem of which the following is  among the most interesting.

\begin{prob}[Nonstandard smooth realization]
Given a diffeomorphism $g$ of a compact manifold $N$  preserving a smooth volume $\mu$ and ergodic with respect to it, does there exist a diffeomorphism  $f$ of a compact manifold $M$
preserving  a smooth volume $\nu$ such that $(f,\nu)$ is measurably isomorphic  to $(g, \mu)$
but not smoothly conjugate to it?
\end{prob}

The only known unconditional restriction to smooth realization
is finiteness of entropy (see \cite[Corollary 3.2.10, Theorem 4.5.3]{KH}).
In other words, no measure preserving transformation with finite entropy is known not to be 
realizable on {\em any} manifold. There are restrictions in low dimension, e.g. any weakly mixing surface diffeomorphism is Bernoulli.

The approximation by conjugation method  provides instances of both nonstandard
 smooth realization  of certain diffeomorphisms and smooth realization of 
maps whose original description is not smooth. This included nonstandard  
smooth realization of  rotations with some Liouvillean rotation numbers
discussed above in \ref{nonstandardrotation} 
(See \cite[Section 4]{AK}; this was recently improved
to any Liouvillean rotation number \cite{SW}).  Very likely  these methods 
should allow  nonstandard realization of 
some Diophantine rotations  by diffeomorphisms of finite differentiability. It looks though that the following problem  tests the limits of power of the approximation by conjugation method and most likely requires  new ideas.

\begin{prob} Does there exist  a nonstandard smooth realization of a Diophantine rotation of the circle by a $C^{\infty}$ diffeomorphism?
\end{prob}

Such a diffeomorphism must act on  a manifold  of dimension  at least two.

In \cite[Section 6]{AK} nonstandard smooth realizations of certain 
translations on tori  were  constructed, as well as smooth realizations of 
certain rotations on  the infinite dimensional torus.
The latter was the earliest example of  smooth realization  of a  
transformation whose natural
``habitat'' is not a finite dimensional manifold. 
In \cite[Section 8]{K1} a general setting is 
described which allows  (with further specification of parameters)  smooth realization of transformations  belonging to measurable combinatorial constructions.
Let us list several characteristic problems which in our view can be 
approached by 
a version of the approximation  by conjugation  method or its modification 
with  decreasing  chances of success.

\begin{prob}
Find a  smooth realization of:
\begin{enumerate}
\item a Gaussian dynamical system with simple (Kronecker) spectrum;
\item a dense $G_{\delta}$ set  of minimal interval exchange transformations;
\item an  adding machine;
\item  the time-one map of  the horocycle flow \ref{horocycleflow} on the modular surface 
$SO(2)\backslash SL(2,\RR)/SL(2,\ZZ)$ (which is not compact, so the 
standard realization cannot be used).
\end{enumerate}
\end{prob}

A version of  the smooth realization problem addresses possibility of  obtaining certain
properties invariant under measurable conjugacy  for volume preserving diffeomorphisms of compact manifolds. In particular one may ask about properties which are generic   for measure preserving transformations in weak  topology  (which are also then generic  for volume preserving homeomorphisms in $C^0$ topology). While  such properties  cannot be generic  in $C^r$ smooth topology by  reasons  quite different  for $r=1$ and $r\ge 2$, one can  consider spaces
of  closures of conjugacies  as in Section \ref{closurespace} and try to determine whether  generic properties of measure preserving diffeomorphisms can be realized there or are generic in those spaces which we denote by $\A$ with various indexes. It seems that the notion of  periodic approximation of  a given type  in the sense of \cite[Definition 1.9]{K1} is quite useful in this respect.  Since approximation of any type with any (arbitrary fast) speed is generic  in weak topology \cite[Theorem 2.1]{K1} any property which follows from any finite or countable combination of approximation type properties is generic too.
Cyclic approximation  where approximating periodic processes have a single tower  is generic in the space $\A$ and certain types of homogeneous approximation   where  there may be many towers  but their heights are all equal \cite[Section 5]{K1} can be shown to be generic too. However periodic processes with towers of unequal height such as type $(n,n+1)$ approximations with 
two towers whose heights differ by one do not naturally appear in the conjugacies of elements of $S^1$ actions. Those approximations and their generalizations are useful is establishing
genericity of various properties (see e.g. \cite[Section 3.3]{K1}).
It looks likely  that in the  case of actions of $\TT^k,\,\, k\ge 2$  it 
may be easier to produce these 
types of approximation in the closure of conjugacies.

\begin{prob}
Given a  circle action $\mS$   and the corresponding space $\A$ is there a diffeomorphism
$f\in\A$ with either of the following properties:
\begin{enumerate}
\item good approximation of type $(n,n+1)$;
\item maximal spectral type disjoint with its convolutions;
\item homogeneous spectrum of multiplicity two for  the Carthesian square $f\times f$?
\end{enumerate}
\end{prob}

See \cite[Sections 3,4]{K1} for relevant  definitions and results.

\begin{prob} Is existence of a measurable square root (or roots of all orders) generic in $\A$?
(Cf. \cite{Ki}).
\end{prob}

\begin{rema}Notice that unlike the previous problem  existence is known here because for example
any rotation of the circle has roots of all orders and there are diffeomorphisms  measurably isomorphic to rotations \cite[section 4]{AK}.
\end{rema}

Given an  action  $\mathcal T$  of  $\TT^2$ denote the closure  of 
conjugacies of elements of this action
in the $C^{\infty}$ topology by $\B$.

\begin{prob}   Is there a diffeomorphism
$f\in\B$ with either of the following properties:
\begin{enumerate}
\item good approximation of type $(n,n+1)$;
\item maximal spectral type disjoint with its convolutions;
\item homogeneous spectrum of multiplicity two for the Carthesian square $f\times f$;
\item mixing?
\end{enumerate}
\end{prob}
Positive answer to the last question will be  of course a discrete time version of Theorem \ref{4 theorem}.  See  Remark \ref{remarkt2mixing}.

\frenchspacing
\bibliographystyle{plain}

\end{document}